\documentclass[10pt]{article}
\usepackage{amsmath,amsthm,amssymb,amsbsy,amsopn}
\usepackage{graphicx}

\newtheorem{Def}{Definition}[section]
\newtheorem{Thm}[Def]{Theorem}
\newtheorem{Lem}[Def]{Lemma}
\newtheorem{Prop}[Def]{Proposition}

\newtheorem{Rem}[Def]{Remark}
\newtheorem{Exam}[Def]{Example}

\theoremstyle{definition}

\newenvironment{enumerate2}{%
   \begin{list}%
   {%
   }%
   {%
      \usecounter{enum2}
      \setlength{\itemindent}{0em}
      \setlength{\leftmargin}{3em}
      \setlength{\rightmargin}{0em}
      \setlength{\labelsep}{1em}
      \setlength{\labelwidth}{3em}
      \setlength{\itemsep}{0em}
      \setlength{\parsep}{0em}
      \setlength{\listparindent}{0em}
   }
}{%
   \end{list}%
}

\makeatletter
\@addtoreset{equation}{section}
\makeatother

\title{On Layered Stable Processes}
\author{C. Houdr\'e\footnote{Laboratoire d'Analyse et de Math\'ematiques
Appliqu\'ees, CNRS UMR 8050, Universit\'e Paris XII, 94010 Cr\'eteil
Cedex, France, and School of Mathematics, Georgia Institute of Technology,
Atlanta, GA, 30332-0160, USA, houdre@math.gatech.edu} 
~and R. Kawai\footnote{Quantitative Research Department, Daiwa Securities SMBC
Co.Ltd., 1-14-5, Eitai, Koto-ku, Tokyo, 135-0034, Japan, reiichiro.kawai@daiwasmbc.co.jp}}
\date{March 28, 2005}

\begin{document}
\maketitle

\begin{abstract}

{\it Layered stable} (multivariate) distributions and processes are
 defined and studied.
A layered stable process combines stable trends of two different
 indices, one of them possibly Gaussian.
More precisely, in short time, it is close to a stable process while, in long
 time, it approximates another stable (possibly Gaussian) process.
We also investigate the absolute continuity of a layered stable process with
 respect to its short time limiting stable process.
A series representation of layered stable processes is derived, giving
 insights into both the structure of the sample paths and of the short
 and long time behaviors.
This series is further used for sample paths simulation.
\end{abstract}


\section{Introduction and preliminaries}

Stable processes form one of the simplest class of L\'evy processes without
Gaussian component.
They have been thoroughly studied by many authors and have
been used in several fields of applications, such as statistical
physics, queueing theory, mathematical finance.
One of their major attractions is the scaling property induced by the
structure of the corresponding L\'evy measure.
Sato~\cite{sato} and Samorodnitsky and Taqqu~\cite{samorodnitsky taqqu}
contain many basic facts on stable distributions and processes.
Recent generalizations of stable processes can also be found, for
example, in Barndorff-Nielsen and Shepard~\cite{barndorff-nielsen
shepard} and in Rosi\'nski~\cite{super rosinski}.
These new classes are also of great interest in applications and have
moreover motivated our study. 

In the present paper, we introduce and study further generalizations which
we call {\it layered stable} distributions and processes.
They are defined in terms of the structure of their L\'evy measure
whose radial component behaves asymptotically as an inverse polynomial
of different orders near zero and at infinity.
{\it The inner and outer (stability) indices} correspond respectively to
these orders of polynomial decay.
This simple layering leads to the following properties:
The outer index determines the moment properties (Proposition
\ref{moment proposition}), while the variational properties depend on
the inner index (Proposition \ref{variation proposition}).
On the other hand, the inner and outer indices also correspond
to short and long time behavior of the sample paths.
In short time, a layered stable process behaves like a stable
process with the corresponding inner index (Theorem \ref{short of eS}).
The long time behavior has two modes depending on the outer index. 
When the outer index is strictly smaller than two, a layered stable
process is close to a stable process with this index, while
behaving like a Brownian motion if the outer index is strictly greater
than two (Theorem \ref{lS long-time theorem}).
In relation to the short time behavior, we investigate the mutual absolute
continuity of a layered stable process and of its short time limiting
stable process (Theorem \ref{absolute continuity theorem}).
A shot noise series representation reveals the nature of layering and
also gives direct insights into the properties of layered stable processes.
We present typical sample paths of a layered stable process, which are
simulated via the series representation for various combinations of stability
indices in order to cover all the types of short and long time behavior.

Let us begin with some general notations which will be used throughout
the text.
$\mathbb{R}^d$ is the $d$-dimensional Euclidean space with the norm
$\|\cdot\|$, $\mathbb{R}_0^d:=\mathbb{R}^d\setminus \{0\},$
$\mathcal{B}(\mathbb{R}^d_0)$ is the Borel $\sigma$-field of
$\mathbb{R}^d_0$, and $S^{d-1}:=\{x\in \mathbb{R}^d:\|x\|=1\}$.
$A'$ is the transpose of the matrix $A$, while $\|\cdot\|_{{\rm o}}$ is
the operator norm of the linear transformation $A\in \mathbb{R}^{d\times
d}$, i.e., $\|A\|_{{\rm o}}=\sup_{\|x\|\le 1}\|Ax\|.$
$f(x) \sim g(x)$ indicates that $f(x)/g(x)\to 1$, as $x\to x_0\in
[-\infty,\infty]$, while $f(x)\asymp g(x)$ is used to mean that there
exist two positive constants $c_1$ and $c_2$ such that $c_1g(x)\le
f(x)\le c_2g(x)$, for all $x$ in an approximate set.
$\mathcal{L}(X)$ is the law of the random vector $X$, while
$\stackrel{\mathcal{L}}{=}$ and $\stackrel{\mathcal{L}}{\to}$
denote, respectively, equality and convergence in distribution, or of
the finite dimensional distributions when random processes are considered.
$\stackrel{d}{\to}$ is used for the weak convergence of
random processes in the space $\mathbb{D}([0,\infty),\mathbb{R}^d)$ of
c\`adl\`ag functions from $[0,\infty)$ into $\mathbb{R}^d$ equipped with
the Skorohod topology.
$\stackrel{v}{\to}$ denotes convergence in the vague topology.
For any $r>0$, $T_r$ is a transformation of measures on $\mathbb{R}^d$
given, for any positive measure $\rho$, by $(T_r\rho)(B)=\rho (r^{-1}B),$
 $B\in\mathcal{B}(\mathbb{R}^d).$
$\mathbb{P}|_{\mathcal{F}_t}$ is the restriction of a
probability measure $\mathbb{P}$ to the $\sigma$-field $\mathcal{F}_t,$
while $\Delta X_t$ denotes the jump of $X$ at time $t,$ that is, $\Delta
X_t:=X_t-X_{t-}.$
Finally, and throughout, all the multivariate or matricial integrals
are defined componentwise. 

Recall that an infinitely divisible probability measure $\mu$ on
$\mathbb{R}^d$, without Gaussian component, is called {\it stable} if
its L\'evy measure is given by
\[
\nu (B)=\int_{S^{d-1}}\sigma (d\xi
)\int_0^{\infty}{\bf 1}_B(r\xi)\frac{dr}{r^{1+\alpha}},\quad B\in
{\mathcal B}(\mathbb{R}_0^d),
\]
where $\alpha \in (0,2)$ is the stability index and where $\sigma$ is a
finite positive measure on $S^{d-1}$.
It is well known that the characteristic function of $\mu$ is given by
\begin{eqnarray}
\widehat{\mu}(y)&=&\exp\left[i\langle y,\eta\rangle +\int_{\mathbb{R}_0^d}
(e^{i\langle y,z\rangle}-1-i\langle y,z\rangle {\bf 1}_{\{\|z\|\le 1\}}(z))
\nu_{\alpha}(dz)\right] \label{stable cf}\\
&=&
\begin{cases}
\exp\left[i\langle y,\tau_{\alpha} \rangle-c_{\alpha}\int_{S^{d-1}}|\langle y,\xi\rangle
|^{\alpha}\left(1-i\tan\frac{\pi \alpha}{2}{\rm sgn}\langle y,\xi
\rangle\right)\sigma (d\xi)\right],&{\rm if}~\alpha \ne 1,\\
\exp\left[i\langle y,\tau_1\rangle-c_1\int_{S^{d-1}}\left(|\langle y,\xi\rangle|+i\frac{2}{\pi}\langle y,\xi\rangle \ln |\langle y,\xi\rangle
|\right)\sigma (d\xi)\right],&{\rm if}~\alpha = 1,
\end{cases}\nonumber
\end{eqnarray}
for some $\eta\in \mathbb{R}^d$, and where $c_{\alpha}=|\Gamma
(-\alpha)\cos \frac{\pi \alpha}{2}|$ when $\alpha \ne 1$ while
$c_1=\pi/2$, with moreover $\tau_{\alpha}=\eta
-\frac{1}{1-\alpha}\int_{S^{d-1}}\xi\sigma (d\xi)$ when $\alpha \ne 1$
and $\tau_1=\eta -(1-\gamma)\int_{S^{d-1}}\xi \sigma (d\xi)$, $\gamma
(=0.5772...)$ being the Euler constant.
A L\'evy process $\{X_t:t\ge 0\}$ such that $\mathcal{L}(X_1)\sim \mu$
is called a {\it stable process}.
Stable processes enjoy the {\it selfsimilarity} property, i.e.,
for any $a>0$,
\[
\{X_{at}:t\ge 0\}\stackrel{\mathcal{L}}{=}\{a^{1/\alpha}X_t+bt:t\ge 0\},
\]
for some $b\in \mathbb{R}^d$.
Next, we recall a shot noise series representation of stable
processes on a fixed finite horizon $[0,T]$, $T>0$. 
Related results can be found, for example, in Theorem 1.4.5 of
Samorodnitsky and Taqqu~\cite{samorodnitsky taqqu}.
The centering constants given below are obtained in Proposition 5.5 of
Rosi\'nski~\cite{super rosinski}.

\begin{Lem}\label{alpha series}
Let $T>0$.
Let $\{T_i\}_{i\ge 1}$ be a sequence of iid uniform random variables on
 $[0,T]$, let $\{\Gamma_i\}_{i\ge 1}$ be an arrival times of a standard Poisson
 process, and let $\{V_i\}_{i\ge 1}$ a sequence of iid random vectors in
 $S^{d-1}$ with common distribution $\sigma (d\xi)/\sigma (S^{d-1})$.
Also let
\begin{equation*}
z_0=
\begin{cases}
0,&if~\alpha \in (0,1),\\
\int_{S^{d-1}}\xi\sigma (d\xi)/\sigma (S^{d-1}) ,&if~\alpha \in [1,2),
\end{cases}
\end{equation*}
and
\begin{equation*}
 b_T=
\begin{cases}
0,&if~\alpha \in (0,1),\\
\sigma (S^{d-1})T(\gamma +\ln (\sigma (S^{d-1})T)),&if~\alpha =1,\\
\left(\frac{\alpha}{\sigma(S^{d-1})T}\right)^{-1/\alpha}\zeta (1/\alpha),&if~\alpha \in (1,2),
\end{cases}
\end{equation*}
where $\zeta$ denotes the Riemann zeta function.
Then, the stochastic process
\[
\left\{\sum_{i=1}^{\infty}\left[\left(\frac{\alpha
	\Gamma_i}{\sigma (S^{d-1})T}\right)^{-1/\alpha}V_i\,{\bf 1}(T_i\le
  t)-\left(\frac{\alpha i}{\sigma (S^{d-1})T}\right)^{-1/\alpha}
z_0\frac{t}{T}\right]+b_Tz_0\frac{t}{T}:t\in [0,T]\right\},
\]
converges almost surely uniformly in $t$ to an $\alpha$-stable process
 $\{X_t:t\in [0,T]\}$ satisfying $\mathbb{E}[e^{i\langle
 y,X_T\rangle}]=\widehat{\mu}(y)^T$, where $\widehat{\mu}$ given by
 (\ref{stable cf}) with
\[
 \eta=
\begin{cases}
\frac{1}{1-\alpha}\int_{S^{d-1}}\xi \sigma (d\xi),&{\rm if}
 ~\alpha \ne 1,\\
0,&{\rm if} ~\alpha =1.
\end{cases}
\]
\end{Lem}

\section{Definition and basic properties}

We first define a layered stable multivariate distribution by precising
the structure of its L\'evy measure in polar coordinates.

\begin{Def}
Let $\mu$ be an infinitely divisible probability measure on
 $\mathbb{R}^d$ and without Gaussian component.
Then, $\mu$ is called layered stable if its L\'evy measure on
 $\mathbb{R}_0^d$ is given by
\begin{equation}
\nu(B)=\int_{S^{d-1}}\sigma
 (d\xi)\int_0^{\infty}{\bf 1}_B(r\xi)q(r,\xi)dr,\quad B\in 
\mathcal{B}(\mathbb{R}_0^d),\label{def of eS}
\end{equation} 
where $\sigma$ is a finite positive measure on $S^{d-1}$, and
 $q$ is a measurable function from $(0,\infty)\times S^{d-1}$ to
 $(0,\infty)$ such that for each $\xi\in S^{d-1}$, 
\begin{equation}\label{short q}
 q(r,\xi)\sim c_1(\xi)r^{-\alpha -1},\quad as~r\to 0,
\end{equation}
and
\begin{equation}\label{long q}
 q(r,\xi)\sim c_2(\xi)r^{-\beta -1},\quad as~r\to\infty,
\end{equation}
where $c_1$ and $c_2$ are integrable (with respect to $\sigma$)
 functions on $S^{d-1}$, and where $(\alpha ,\beta)\in (0,2)\times
 (0,\infty )$.
\end{Def}

$q(\cdot ,\cdot )$ is called {\it the $q$-function} of $\mu$, or of its
L\'evy measure $\nu$.  
Clearly, $\nu$ is well defined as a
L\'evy measure since it behaves like an $\alpha$-stable L\'evy measure
near the origin while decaying like a $\beta$-Pareto density when
sufficiently far away from the origin.
$\alpha$ and $\beta$ are respectively called {\it the inner and
outer (stability) indices} of $\mu$, or of $\nu$.

{\it For convenience, we henceforth use the notations $\sigma_1$
and $\sigma_2$ for the finite positive measures on $S^{d-1}$ defined
respectively by 
\begin{equation}\label{def of sigma1}
\sigma_1(B):=\int_Bc_1(\xi)\sigma (d\xi),\quad B\in \mathcal{B}(S^{d-1}), 
\end{equation}
and
\begin{equation}\label{def of sigma2}
\sigma_2(B):=\int_Bc_2(\xi)\sigma (d\xi),\quad B\in \mathcal{B}(S^{d-1}), 
\end{equation}
while $\nu^{\alpha}_{\sigma}$ is used for the positive measure on
$\mathbb{R}_0^d$ given by
\begin{equation}\label{def of stable levy measure}
 \nu^{\alpha}_{\sigma}(B):=\int_{S^{d-1}}\sigma (d\xi)\int_0^{\infty}{\bf
 1}_B(r\xi)\frac{dr}{r^{\alpha+1}},\quad B\in \mathcal{B}(\mathbb{R}_0^d),
\end{equation}
where $\alpha\in (0,\infty)$ and where $\sigma$ is a finite positive 
measure on $S^{d-1}$.}
Note that if $\alpha \in (0,2)$, $\nu_{\sigma}^{\alpha}$ is simply an
$\alpha$-stable L\'evy measure, while not well defined as a L\'evy
measure when $\alpha \ge 2$.

\begin{Exam}\label{special remark}{\rm 
The following layered stable L\'evy measure is simple, yet interesting:
\begin{eqnarray}
 \nu (B)&=&\int_B{\bf 1}_{\{\|z\|\le
  1\}}(z)\nu_{\sigma}^{\alpha}(dz)+\int_B{\bf 1}_{\{\|z\|>
  1\}}(z)\nu_{\sigma}^{\beta}(dz)\nonumber \\
&=&\int_{S^{d-1}}\sigma (d\xi)\int_0^{\infty}{\bf
 1}_B(r\xi)\frac{dr}{r^{\alpha+1}{\bf 1}_{(0,1]}(r)+r^{\beta+1}{\bf
 1}_{(1,\infty)}(r)},\quad B\in
 \mathcal{B}(\mathbb{R}_0^d).\label{original lS}
\end{eqnarray}
The corresponding $q$-function is given by
\[
 q(r,\xi)=\sigma(S^{d-1})^{-1}(r^{-\alpha-1}{\bf 1}_{(0,1]}(r)+r^{-\beta-1}{\bf
 1}_{(1,\infty)}(r)),\quad \xi\in S^{d-1},
\]
which is independent of $\xi$.
The measure $\nu$ consists of two disjoint domains of stability,
and this construction results in two layers for the radial
 component associated with each respective stability index.
The name ``{\it layered stable}'' originates from this special
 structure.

Recall that an infinitely divisible probability measure $\mu$ on
 $\mathbb{R}^d$ is said to be of class $L_0$, or
 selfdecomposable if for any $b>1$, there exists a probability measure
 $\varrho_b$ such that
 $\widehat{\mu}(z)=\widehat{\mu}(b^{-1}z)\widehat{\varrho_b}(z).$
Equivalently, the L\'evy measure of $\mu$ has the form
\[
\int_{S^{d-1}}\sigma(d\xi)\int_0^{\infty}{\bf
 1}_B(r\xi)k_{\xi}(r)\frac{dr}{r},\quad B\in \mathcal{B}(\mathbb{R}_0^d),
\]
where $\sigma$ is a finite positive measure on $S^{d-1}$ and where
 $k_{\xi}(r)$ is a nonnegative function measurable in $\xi \in S^{d-1}$ and
 decreasing in $r>0.$
Clearly, the L\'evy measure (\ref{original lS}) induces a
 selfdecomposable measure.
Moreover, the classes $L_m$, $m=1,2,\ldots ,$ are defined recursively as
 follows; $\mu\in L_m$ if for every $b>1$, there exists $\varrho_b\in
 L_{m-1}$ such that 
 $\widehat{\mu}(z)=\widehat{\mu}(b^{-1}z)\widehat{\varrho_b}(z).$
Clearly, $L_0\supset L_1 \supset L_2 \supset \cdots$.
Let $h_{\xi}(u):=k_{\xi}(e^{-u})$, be the so-called {\it $h$-function} of
 $\mu$, or of its L\'evy measure.
Then, alternatively, $\mu \in L_0$ is shown to be in $L_m$ if and only if
 $h_{\xi}(u)\in C^{m-1}$ and $h^{(j)}\ge 0$, for $j=0,1,\ldots ,m-1$.
(See Sato \cite{sato class L} for more details.)
The $h$-function of the L\'evy measure (\ref{original lS}) is given by
\[
 h_{\xi}(u)=e^{\alpha u}{\bf 1}_{(0,\infty)}(u)+e^{\beta u}{\bf 1}_{(-\infty,0]}(u),
\]
which is in $C^0$ but not in $C^1$.
Therefore, the infinitely divisible probability measure induced by
 (\ref{original lS}) is in $L_1$, but not in $L_2$.
}\end{Exam}

The following result asserts that a layered stable distribution has the
same probability tail behavior as $\beta$-Pareto distributions, or
$\beta$-stable distributions if $\beta\in (0,2)$.
 
\begin{Prop}\label{moment proposition}{\rm (Moments)}
Let $\mu$ be a layered stable distribution with L\'evy measure $\nu$
 given by (\ref{def of eS}) and let $\sigma_2$ be the measure (\ref{def
 of sigma2}).
If $\sigma_2(S^{d-1})\ne 0$, then
\[
 \int_{\mathbb{R}^d}\|x\|^{p}\mu (dx)
\begin{cases}
<\infty, & p\in (0,\beta),\\
=\infty, & p\in [\beta,\infty).
\end{cases}
\]
Moreover, $\int_{\mathbb{R}^d}\|x\|^p\mu (dx)<\infty$, $p\ge \beta$ and
 $\int_{\mathbb{R}^d}e^{\theta \|x\|}\mu (dx)<\infty$, $\theta >0$ if
 and only if $\sigma_2(S^{d-1})=0.$
\end{Prop}

\begin{proof}
By Theorem 25.3 of Sato \cite{sato}, it is enough to show that the
 restriction of $\nu$ to the set $\{z\in \mathbb{R}_0^d:\|z\|>1\}$ has
 the corresponding moment properties.

First, assume $\sigma_2(S^{d-1})\ne 0$.
Observe that $\int_{\|z\|>1}\|z\|^p \nu(dz)=
 \int_{S^{d-1}}\sigma(d\xi)\int_1^{\infty}r^pq(r,\xi)dr$, and then by
 (\ref{long q}), the right hand side is bounded from above and below by
 constant multiples of
 $\sigma_2(S^{d-1})\int_1^{\infty}r^p\frac{dr}{r^{\beta +1}}$ if $p\in
 (0,\beta)$, while it is otherwise clearly infinite.

Next, assume $\sigma_2(S^{d-1})=0$ and let $p\in [\beta,\infty)$.
Then, there exists $M>0$ such that
 $\int_{\|z\|>1}\|z\|^p\nu (dz)\asymp\int_{S^{d-1}}\sigma(d\xi)\int_1^Mr^pq(r,\xi)dr$ and
 $\int_{\|z\|>1}e^{\theta \|z\|}\nu (dz)\asymp
 \int_{S^{d-1}}\sigma(d\xi)\int_1^Me^{\theta r}q(r,\xi)dr.$
Conversely, if $\sigma_2(S^{d-1})\ne 0$ and $p\in [\beta,\infty)$,
then $\int_{\|z\|>1}\|z\|^p\nu (dz)=+\infty$ as already shown and, again
 by (\ref{long q}), $\int_{\|z\|>1}e^{\theta \|z\|}\nu
 (dz)=\int_{S^{d-1}}\sigma (d\xi)\int_1^{\infty}e^{\theta
 r}q(r,\xi)dr=+\infty.$
\end{proof}

Let us define the associated L\'evy processes.

\begin{Def}
A L\'evy process, without Gaussian component, is called layered stable
 if its L\'evy measure is given by (\ref{def of eS}).
\end{Def}

Henceforth, $\{X^{LS}_t:t\ge 0\}$ denotes a layered stable process
in $\mathbb{R}^d$.
Its characteristic function at time $1$ is given by
\begin{equation}\label{def via cf}
 \mathbb{E}[e^{i\langle y,X^{LS}_1\rangle}]=\exp\left[i\langle
 y,\eta\rangle +\int_{\mathbb{R}_0^d}
(e^{i\langle y,z\rangle }-1-i\langle y,z\rangle
 {\bf 1}_{\{\|z\|\le 1\}}(z))\nu(dz)\right],
\end{equation}
where $\nu$ is the L\'evy measure given by (\ref{def of eS}) and
$\eta \in\mathbb{R}^d$.
For convenience of notation, we write $\{X^{LS}_t:t\ge 0\}\sim
LS_{\alpha,\beta}(\sigma,q;\eta)$ when (\ref{def via cf}) holds.
Similarly, for $\alpha\in (0,2)$, $\{X^{(\alpha)}_t:t\ge 0\}$ denotes an
$\alpha$-stable L\'evy process.
Its characteristic function at time $1$ is given by
\begin{equation}\label{cf of stable for lS}
 \mathbb{E}[e^{i\langle y,X^{(\alpha)}_1\rangle}]=
\begin{cases}
\exp\left[i\langle y,\eta\rangle+\int_{\mathbb{R}_0^d}(e^{i\langle
 y,z\rangle}-1)\nu^{\alpha}_{\sigma}(dz)\right],&if~\alpha\in (0,1),\\
\exp\left[i\langle y,\eta\rangle+\int_{\mathbb{R}_0^d}(e^{i\langle
 y,z\rangle}-1-i\langle y,z\rangle
 {\bf 1}_{\{\|z\|\le 1\}}(z))\nu^1_{\sigma}(dz)\right],&if~\alpha =1,\\
\exp\left[i\langle y,\eta\rangle+\int_{\mathbb{R}_0^d}(e^{i\langle
 y,z\rangle}-1-i\langle y,z\rangle
 )\nu^{\alpha}_{\sigma}(dz)\right],&if~\alpha\in (1,2),
\end{cases}
\end{equation}
where $\nu_{\sigma}^{\alpha}$ is given by (\ref{def of stable levy measure}),
and we write $\{X^{(\alpha)}_t:t\ge 0\}\sim S_{\alpha}(\sigma;\eta)$ when
(\ref{cf of stable for lS}) holds.  

\vspace{1em}
A layered stable process shares the variational
properties of a stable process with inner index $\alpha$.

\begin{Prop}\label{variation proposition}{\rm ($p$-th variation)}
Let $X:=\{X^{LS}_t:t\ge 0\}\sim LS_{\alpha,\beta}(\sigma,q;\eta)$.

\noindent (i) If $\sigma_1(S^{d-1})>0,$ then $X$ is a.s. of finite variation on
 every interval of positive length if and only if $\alpha \in (0,1).$

\noindent (ii) If $\sigma_1(S^{d-1})>0$, $(\alpha,\beta) \in [1,2)\times
 (1,\infty)$ and $\eta =-\int_{S^{d-1}}\xi\sigma
 (d\xi)\int_1^{\infty}rq(r,\xi)dr,$ then $X$ is a.s. of finite $p$-th
 variation on every interval of positive length if and only if $p> \alpha$.

\noindent (iii) If $\sigma_1(S^{d-1})=0,$ then it is a.s. of finite
 variation on every interval of positive length.
\end{Prop}

\begin{proof}
(i) Recall that the radial component of the layered stable L\'evy measure
 near the origin behaves like the one of an $\alpha$-stable L\'evy measure.  
The first claim then follows immediately from Theorem 3 of Gikhman and
 Skorokhod \cite{gikhman skorokhod}.

(ii) Since $X$ is now centered, Th\'eor\`eme III b of Bretagnolle
 \cite{bretagnolle} directly applies.

(iii) Letting $\nu$ be the L\'evy measure of $X$, there exists
 $\epsilon\in (0,1)$ such that $\nu (\{z\in \mathbb{R}_0^d:\|z\|\le
 \epsilon\})< \infty$ and so $\int_{\|z\|\le 1}\|z\|^p\nu (dz)<\infty$,
 $p\ge 1.$
As in (i), the result follows from Theorem 3 of Gikhman and Skorokhod
 \cite{gikhman skorokhod}.
\end{proof}

Let us now consider a series representation for a general layered  
stable process $\{X^{LS}_t:t\ge 0\}\sim LS_{\alpha,\beta}(\sigma,q;0).$
Fix $T>0$.
Let $\{T_i\}_{i\ge 1}$ be a sequence of iid
 uniform random variables on $[0,T]$, let $\{\Gamma_i\}_{i\ge 1}$ be
 Poisson arrivals with rate $1$, and let $\{V_i\}_{i\ge 1}$ be a
 sequence of iid random vectors in $S^{d-1}$ with common
 distribution $\sigma (d\xi)/\sigma (S^{d-1})$. 
Assume moreover that the random sequences $\{T_i\}_{i\ge 1}$,
$\{\Gamma_i\}_{i\ge 1}$, and $\{V_i\}_{i\ge 1}$ are all mutually
independent.
Also, let
\[
  \overleftarrow{q}(u,\xi):=\inf\{r>0:q([r,\infty),\xi)<u\},
\]
and let $\{b_i\}_{i\ge 1}$ be a sequence of constants given by 
\[
 b_i=\int_{i-1}^i\mathbb{E}[\overleftarrow{q}(s/T,V_1)V_1{\bf
 1}(\overleftarrow{q}(s/T,V_1)\le 1)]ds,
\]
Then, by Theorem 5.1 of Rosi\'nski \cite{rosinski2} with the help of
the LePage's method \cite{lepage}, the stochastic process
\begin{equation}\label{original series}
 \left\{\sum_{i=1}^{\infty}\left[\overleftarrow{q}(\Gamma_i/T,V_i)V_i{\bf
 1}(T_i\le t)-b_i\frac{t}{T}\right]:t\in [0,T]\right\},
\end{equation}
converges almost surely uniformly in $t$ to a L\'evy process whose
 marginal law at time $1$ is $LS_{\alpha,\beta}(\sigma,q;0).$

\begin{Exam}{\rm 
The L\'evy measure (\ref{original lS}) leads to a very illustrative
 series representation.
Indeed,
\[
 \overleftarrow{q}(r,\xi)=\left(\frac{\beta r}{\sigma
 (S^{d-1})}\right)^{-1/\beta}{\bf
 1}_{(0,\sigma(S^{d-1})/\beta]}(r)+\left(\frac{\alpha r}{\sigma
 (S^{d-1})}+1-\frac{\alpha}{\beta}\right)^{-1/\alpha}{\bf 1}_{(\sigma
 (S^{d-1})/\beta,\infty)}(r),
\]
and so the stochastic process
\begin{eqnarray}
&&\Bigg\{\sum_{i=1}^{\infty}
\Bigg[\Bigg(\left(\frac{\beta\Gamma_i}{\sigma
(S^{d-1})T}\right)^{-1/\beta}{\bf 1}_{(0,\sigma
(S^{d-1})T/\beta]}(\Gamma_i)\label{special series}\\
&&~+\left(\frac{\alpha\Gamma_i}{\sigma
   (S^{d-1})T}+1-\frac{\alpha}{\beta}\right)^{-1/\alpha}{\bf 1}_{(\sigma
 (S^{d-1})T/\beta,\infty)}(\Gamma_i)\Bigg)V_i{\bf 1}(T_i\le t)
-b_iz_0\frac{t}{T}\Bigg]:t\in [0,T]\Bigg\},\nonumber
\end{eqnarray}
where
\[
b_i=\left(\frac{\beta}{\sigma(S^{d-1})T}\right)^{-1/\beta}\frac{(i\land
 \sigma (S^{d-1})T/\beta)^{1-1/\beta}-((i-1)\land \sigma (S^{d-1})T/\beta)^{1-1/\beta}}{1-1/\beta},
\]
converges almost surely uniformly in $t$ to a L\'evy process whose
 marginal law at time $1$ is $LS_{\alpha,\beta}(\sigma,q;0)$, with
 $z_0=\int_{S^{d-1}}\xi\sigma(d\xi)/\sigma(S^{d-1})$.
This series representation directly reveals the nature of layering; all
 jumps with absolute size greater than $1$ are due to the $\beta$-stable
 shot noise series $\left(\frac{\beta \Gamma_i}{\sigma
 (S^{d-1})}\right)^{-1/\beta}V_i$, while smaller jumps come from
 $\left(\frac{\alpha \Gamma_i}{\sigma
 (S^{d-1})}+1-\frac{\alpha}{\beta}\right)^{-1/\alpha}V_i$, which
 resembles $\alpha$-stable jumps.
}\end{Exam}

\section{Short and long time behavior}
We now present one of the two main results of this
section by giving the short time behavior of a layered stable process.
The results of this section were motivated by Section 3 of Rosi\'nski
\cite{super rosinski}, where stable behavior is obtained (for tempered stable processes) 
in short time while Gaussian convergence is obtained in long time.  
Here, in addition, we also obtain a further level of stable (non--Gaussian) 
convergence in long time.  
Recall that $\sigma_1$ and $\sigma_2$ are the finite positive measures
respectively given in (\ref{def of sigma1}) and (\ref{def of sigma2}),
and that for any $r>0,$ $T_r$ transforms the positive measure $\rho$, via
$(T_r\rho)(B)=\rho (r^{-1}B),$ $B\in \mathcal{B}(\mathbb{R}^d).$ 
For convenience, we will use the notation
$\nu_{\sigma,q}^{\alpha,\beta}$ for the L\'evy measure of 
a layered stable process $LS_{\alpha,\beta}(\sigma,q;\eta)$ throughout
this section.

\begin{Thm}\label{short of eS} \underline{\it Short time behavior}:
Let $\{X^{LS}_t:t\ge 0\}\sim LS_{\alpha,\beta}(\sigma,q;0)$, let
\[
 \eta_{\alpha,\beta}=
\begin{cases}
\int_{S^{d-1}}\xi\sigma (d\xi)\int_0^1r q(r,\xi)dr,&if~\alpha\in (0,1),\\
-\int_{S^{d-1}}\xi\sigma (d\xi)\int_1^{\infty}r
 q(r,\xi)dr,&if~(\alpha,\beta)\in (1,2)\times(1,\infty),\\
0,&otherwise,
\end{cases}
\]
and let
\[
 b_{\alpha,\beta}=
\begin{cases}
\frac{1}{\alpha -1}\int_{S^{d-1}}\xi\sigma_1 (d\xi),&if~(\alpha,\beta) \in
 (1,2)\times (0,1],\\
0,&otherwise.
\end{cases}
\]
Then,
\[
 \{h^{-1/\alpha}(X^{LS}_{ht}+ht\eta_{\alpha,\beta})-t b_{\alpha,\beta}
:t\ge 0\}\stackrel{d}{\to}\{X^{(\alpha)}_t:t\ge 0\},\quad as ~h\to 0,
\]
where $\{X^{(\alpha)}_t:t\ge 0\}\sim S_{\alpha}(\sigma_1;0)$.
\end{Thm}

\begin{proof}
Since a layered stable process is a L\'evy process, by a theorem of
 Skorohod (see Theorem 15.17 of Kallenberg \cite{kallenberg}), it
 suffices to show the weak convergence of its marginals at time $1$.
To this end, we will show the proper convergence of the generating
 triplet of the infinitely divisible law, following Theorem 15.14 of
 Kallenberg \cite{kallenberg}.

For the convergence of the L\'evy measure, we need to show that as $h\to
 0$,
\[
h(T_{h^{-1/\alpha}}\nu^{\alpha,\beta}_{\sigma,q})\stackrel{v}{\to}
\nu^{\alpha}_{\sigma_1}, 
\]
or equivalently,
\[
 \lim_{h\to 0}\int_{\mathbb{R}_0^d}f(z)h(T_{h^{-1/\alpha}}\nu^{\alpha
 ,\beta}_{\sigma,q})(dz)=\int_{\mathbb{R}_0^d}f(z)\nu^{\alpha}_{\sigma_1}(dz),
\]
for all bounded continuous function $f:\mathbb{R}^d_0\to \mathbb{R}$ 
vanishing in a neighborhood of the origin.
Letting $f$ be such a function with $|f|\le C<\infty$ and $f(z)\equiv 0$
 on $\{z\in \mathbb{R}^d_0:\|z\|\le \epsilon\}$, for some $\epsilon>0$,
we get by (\ref{short q}),
\begin{eqnarray*}
\int_{\mathbb{R}_0^d}f(z)h(T_{h^{-1/\alpha}}\nu^{\alpha,\beta}_{\sigma,q})
(dz)&=&\int_{S^{d-1}}\sigma (d\xi)\int_0^{\infty}f(h^{-1/\alpha}r\xi)h
q(r,\xi)dr\\
&=&\int_{S^{d-1}}\sigma (d\xi)\int_0^{\infty}f(r\xi)h^{1+1/\alpha}
q(h^{1/\alpha}r,\xi)dr\\
&\to&\int_{S^{d-1}}c_1(\xi)\sigma(d\xi)\int_0^{\infty}f(r\xi)\frac{dr}{r^{\alpha+1}},
\end{eqnarray*}
as $h\to 0$, where the last convergence holds true since for $h\in (0,1)$,
\[
 \left|\int_{S^{d-1}}\sigma (d\xi)\int_0^{\infty}f(h^{-1/\alpha}r\xi)h
q(r,\xi)dr\right|\le C \left|\int_{S^{d-1}}\sigma (d\xi)\int_{\epsilon}^{\infty}q(r,\xi)dr\right|<\infty.
\]

For the convergence of the Gaussian component, we need to show that for
 each $\kappa >0,$
\[
\int_{\|z\|\le \kappa}zz'h(T_{h^{-1/\alpha}}\nu^{\alpha,\beta}_{\sigma,q})(dz)
\to \int_{\|z\|\le \kappa}zz'\nu^{\alpha}_{\sigma_1}(dz),
\]
as $h\to 0$.
Again, by (\ref{short q}),
\begin{eqnarray*}
 \int_{\|z\|\le
 \kappa}zz'h(T_{h^{-1/\alpha}}\nu^{\alpha,\beta}_{\sigma,q})(dz)
&=&\int_{S^{d-1}}\xi\xi'\sigma
 (d\xi)\int_0^{h^{1/\alpha}\kappa}r^2h^{1-2/\alpha}q(r,\xi)dr\\
&=&\int_{S^{d-1}}\xi\xi'\sigma(d\xi)\int_0^{\kappa}r^2h^{1+1/\alpha}
q(h^{1/\alpha}r,\xi)dr\\
&\to&
 \int_{S^{d-1}}\xi\xi'\sigma_1(d\xi)\int_0^{\kappa}r^2\frac{dr}{r^{\alpha+1}}\\
&=&\int_{\|z\|\le\kappa}zz'\nu^{\alpha}_{\sigma_1}(dz),
\end{eqnarray*}
where the passage to the limit is justified since, for $h\in (0,1),$
\[
 \left\|\int_{S^{d-1}}\xi\xi'\sigma (d\xi)\int_0^{h^{1/\alpha}\kappa}r^2
h^{1-2/\alpha}q(r,\xi)dr\right\|_{{\rm o}}\le \left\|\int_{S^{d-1}}\xi\xi'
\sigma (d\xi)\int_0^{\kappa}r^2q(r,\xi)dr\right\|_{{\rm o}} <\infty.
\]

For the convergence of the drift part, assume first that
 $(\alpha,\beta) \notin (1,2)\times (0,1].$ 
For a $\sigma$-finite positive measure $\nu$ on $\mathbb{R}_0^d$, let
\begin{equation}\label{def of constant c}
 C_{\alpha}(\nu):=
\begin{cases}
\int_{\|z\|\le 1}z\nu (dz),&{\rm if}~\alpha \in (0,1),\\
0,&{\rm if}~\alpha =1,\\
-\int_{\|z\|>1}z\nu (dz),&{\rm if}~\alpha \in (1,2).
\end{cases}
\end{equation}
Clearly, $\eta_{\alpha,\beta}=C_{\alpha}(\nu_{\sigma,q}^{\alpha,\beta})$
 and we then show that as $h\to 0$,
\[
C_{\alpha}(h(T_{h^{-1/\alpha}}\nu^{\alpha,\beta}_{\sigma,q}))-
\int_{\kappa <\|z\|\le 1}zh(T_{h^{-1/\alpha}}\nu^{\alpha,
\beta}_{\sigma,q})(dz) \to C_{\alpha}(\nu^{\alpha}_{\sigma_1})-
\int_{\kappa <\|z\|\le 1}z\nu^{\alpha}_{\sigma_1}(dz),
\]
for each $\kappa >0.$
Letting 
\[
 B=
\begin{cases}
\{z\in \mathbb{R}_0^d:\|z\|\le \kappa\},&{\rm if}~\alpha \in (0,1),\\
\{z\in \mathbb{R}_0^d:\kappa <\|z\|\le 1\},&{\rm if}~\alpha =1,\\
\{z\in \mathbb{R}_0^d:\|z\|>\kappa\},&{\rm if}~(\alpha,\beta)\in
 (1,2)\times (1,\infty),\\
\end{cases}
\]
we have as $h\to 0$,
\begin{eqnarray*}
\int_{\mathbb{R}_0^d}{\bf 1}_B(z)zh(T_{h^{-1/\alpha}}\nu_{\sigma,q}^{\alpha,\beta})(dz)
&=&\int_{S^{d-1}}\xi\sigma (d\xi)\int_0^{\infty}{\bf
1}_B(h^{-1/\alpha}r\xi)rh^{1-1/\alpha} q(r,\xi)dr\\
&=&\int_{S^{d-1}}\xi\sigma (d\xi)\int_0^{\infty}{\bf
 1}_B(r\xi)rh^{1+1/\alpha}q(h^{1/\alpha}r,\xi)dr\\
&\to&\int_{S^{d-1}}\xi\sigma_1(d\xi)\int_0^{\infty}{\bf 1}_B(r\xi)r
\frac{dr}{r^{\alpha+1}},
\end{eqnarray*}
where the convergence holds true since for $h\in (0,1)$, and with the
 help of (\ref{short q}),
\[
 \left\|\int_{S^{d-1}}\xi\sigma (d\xi)\int_0^{\infty}{\bf
 1}_B(r\xi)rh^{1+1/\alpha}q(h^{1/\alpha}r,\xi)dr\right\|\asymp
 \left\|\int_{S^{d-1}}\xi\sigma (d\xi)\int_0^{\infty}{\bf
 1}_B(r\xi)rq(r,\xi)dr\right\|<\infty.
\]
Finally, assume $(\alpha,\beta) \in (1,2)\times (0,1].$
Then, as $h\to 0$,
\[
 -b_{\alpha,\beta}-\int_{\kappa<\|z\|\le 1}zh(T_{h^{-1/\alpha}}\nu_{\sigma,q}^{\alpha,\beta})(dz)\to-\int_{\|z\|>\kappa}z\nu_{\sigma_1}^{\alpha}(dz),
\]
for each $\kappa >0$, where the convergence holds true as before.
This completes the proof.
\end{proof}

Our next result is also important.
Unlike in short time, the long time behavior of a layered stable
process depends on its outer stability index $\beta$.
This behavior is akin to a $\beta$-stable process if $\beta \in (0,2)$,
while akin to a Brownian motion whenever $\beta \in (2,\infty ).$

\begin{Thm}\label{lS long-time theorem} \underline{\it Long time
 behavior}: Let $\{X^{LS}_t:t\ge 0\}\sim
 LS_{\alpha,\beta}(\sigma ,q;0).$

\noindent (i) Let $\beta \in (0,2),$ let
\[
 \eta_{\alpha,\beta}=
\begin{cases}
\int_{S^{d-1}}\xi\sigma (d\xi)\int_0^1rq(r,\xi)dr,&if~(\alpha,\beta)\in
 (0,1)\times (0,1),\\
-\int_{S^{d-1}}\xi\sigma (d\xi)\int_1^{\infty}rq(r,\xi)dr,&if~\beta\in
 (1,2),\\
0,&otherwise,
\end{cases}
\]
and let
\[
 b_{\alpha,\beta}=
\begin{cases}
\frac{1}{1-\beta}\int_{S^{d-1}}\xi\sigma_2(d\xi),&if~(\alpha,\beta)\in
 [1,2)\times (0,1),\\
0,&otherwise.
\end{cases}
\]
Then,
\[
 \{h^{-1/\beta}(X^{LS}_{ht}+ht\eta_{\alpha,\beta})+tb_{\alpha,\beta}:t\ge 0\}\stackrel{d}{\to}
\{X^{(\beta)}_t:t\ge 0\},\quad as ~h\to \infty ,
\]
where $\{X^{(\beta)}_t:t\ge 0\}\sim S_{\beta}(\sigma_2;0).$

\noindent (ii) Let $\beta \in (2,\infty )$ and let
\begin{equation}\label{def of constant a}
 \eta=-\int_{S^{d-1}}\xi\sigma(d\xi)\int_1^{\infty}rq(r,\xi)dr.
\end{equation}
Then,
\begin{equation}
 \{h^{-1/2}(X^{LS}_{ht}+ht\eta):t\ge 0\}\stackrel{d}{\to}
\{W_t:t\ge 0\},\quad as ~h\to \infty ,\label{eS long 2}
\end{equation}
where $\{W_t:t\ge 0\}$ is a centered Brownian motion with covariance
 matrix $\int_{\mathbb{R}_0^d}zz'\nu_{\sigma,q}^{\alpha,\beta}(dz).$
\end{Thm} 

\begin{proof}
The claim (i) can be proved as (i) in Theorem \ref{short of eS}.
For the convergence of the L\'evy measure, we will show that
\[
\lim_{h\to
 \infty}\int_{\mathbb{R}_0^d}f(z)h(T_{h^{-1/\beta}}\nu^{\alpha
 ,\beta}_{\sigma,q})(dz)=\int_{\mathbb{R}_0^d}f(z)\nu^{\beta}_{\sigma_2}(dz),
\]
for all bounded continuous function $f:\mathbb{R}^d_0\to
\mathbb{R}$ vanishing in a neighborhood of the origin.
Letting $f$ be such a function with $|f|\le C<\infty$ and $f(z)\equiv 0$ on
 $\{z\in \mathbb{R}^d_0:\|z\|\le \epsilon\}$, for some $\epsilon >0$, we
 get by (\ref{long q}),
\begin{eqnarray*}
\int_{\mathbb{R}_0^d}f(z)h(T_{h^{-1/\beta}}\nu^{\alpha,\beta}_{\sigma,q})
(dz)&=&\int_{S^{d-1}}\sigma (d\xi)\int_0^{\infty}f(h^{-1/\beta}r\xi)h
q(r,\xi)dr\\
&=&\int_{S^{d-1}}\sigma (d\xi)\int_0^{\infty}f(r\xi)h^{1+1/\beta}
q(h^{1/\beta}r,\xi)dr\\
&\to&\int_{S^{d-1}}c_2(\xi)\sigma(d\xi)\int_0^{\infty}f(r\xi)\frac{dr}{r^{\beta+1}},
\end{eqnarray*}
as $h\to \infty$, where the last convergence holds true because of
 (\ref{long q}) and since for sufficiently large $h>0$,
\begin{eqnarray*}
 \left|\int_{S^{d-1}}\sigma (d\xi)\int_0^{\infty}f(h^{-1/\beta}r\xi)h
q(r,\xi)dr\right|&=&\left|\int_{S^{d-1}}\sigma (d\xi)\int_{h^{1/\beta}
\epsilon}^{\infty}f(h^{-1/\beta}r\xi)hq(r,\xi)dr\right|\\
&\le&hC\int_{S^{d-1}}\sigma(d\xi)\int_{h^{1/\beta}\epsilon}^{\infty}q(r,\xi)dr\\&\asymp&hC
 \sigma_2
 (S^{d-1})\int_{h^{1/\beta}\epsilon}^{\infty}\frac{dr}{r^{\beta+1}}\\
&=&C\sigma_2(S^{d-1})\frac{\epsilon^{-\beta}}{\beta}<\infty.
\end{eqnarray*}
For the convergence of the Gaussian component, we have as $h\to \infty$
 and for each $\kappa >0,$
\begin{eqnarray*}
\int_{\|z\|\le
 \kappa}zz'h(T_{h^{-1/\beta}}\nu_{\sigma,q}^{\alpha,\beta})(dz)&=&\int_{S^{d-1}}\xi\xi'\sigma
 (d\xi)\int_0^{\kappa}r^2h^{1+1/\beta}q(h^{1/\beta}r,\xi)dr\\
&\to&\int_{S^{d-1}}\xi\xi'\sigma_2
 (d\xi)\int_0^{\kappa}r^2\frac{dr}{r^{\beta+1}}\\
&=&\int_{\|z\|\le \kappa}zz'\nu_{\sigma_2}^{\beta}(dz),
\end{eqnarray*}
where the passage to the limit is justified next.
Let $h\in (\kappa^{-\beta},+\infty)$ and write
\begin{eqnarray*}
\left\|\int_{\|z\|\le \kappa}zz'h(T_{h^{-1/\beta}}\nu_{\sigma,q}^{
\alpha,\beta})(dz)\right\|_{{\rm o}}&\le &
\left\|\int_{S^{d-1}}\xi\xi'\sigma(d\xi)\int_0^1h^{1-2/\beta}
r^2q(r,\xi)dr\right\|_{{\rm o}}\\
&&+\left\|\int_{S^{d-1}}\xi\xi'\sigma(d\xi)\int_1^{h^{1/\beta}\kappa}h^{1-2/\beta}
r^2q(r,\xi)dr\right\|_{{\rm o}}.
\end{eqnarray*}
The first term of the right hand side above is clearly bounded by
 $\kappa^{2-\beta} \|\int_{S^{d-1}}\xi\xi'\sigma(d\xi)\int_0^1
 r^2q(r,\xi)dr\|_{{\rm o}}$, while the second term is also bounded since
 for $h\in (\kappa^{-\beta},+\infty),$ 
\begin{eqnarray*}
\left\|\int_{S^{d-1}}\xi\xi'\sigma(d\xi)\int_1^{h^{1/\beta}\kappa}h^{1-2/\beta}
r^2q(r,\xi)dr\right\|_{{\rm
o}}&\asymp&h^{1-2/\beta}\int_1^{h^{1/\beta}\kappa}r^2\frac{dr}{r^{\beta+1}}
\left\|\int_{S^{d-1}}\xi\xi'\sigma_2(d\xi)\right\|_{{\rm o}}\\
&=&\frac{\kappa^{2-\beta}-h^{1-2/\beta}}{2-\beta}\left\|\int_{S^{d-1}}
\xi\xi'\sigma_2(d\xi)\right\|_{{\rm o}}<\infty.
\end{eqnarray*}

Finally, we study the convergence of the drift part.
Assume first that $(\alpha,\beta)\notin [1,2)\times (0,1).$
Let $C_{\beta}(\nu)$ be the constant defined as in (\ref{def of constant
 c}) but depending on $\beta$ and $\nu$.
Clearly, $\eta_{\alpha,\beta}=C_{\beta}(\nu_{\sigma,q}^{\alpha,\beta}).$
We will then show that as $h\to \infty$,
\[
C_{\beta}(h(T_{h^{-1/\beta}}\nu^{\alpha,\beta}_{\sigma,q}))-\int_{\kappa
 <\|z\|\le 1}zh(T_{h^{-1/\beta}}\nu^{\alpha,\beta}_{\sigma,q})(dz)\to 
C_{\beta}(\nu^{\beta}_{\sigma_2})-\int_{\kappa <\|z\|\le
 1}z\nu^{\beta}_{\sigma_2}(dz),
\]
for each $\kappa >0$.
Letting 
\[
 B=
\begin{cases}
\{z\in \mathbb{R}_0^d:\|z\|\le \kappa\},&{\rm if}~\beta \in (0,1),\\
\{z\in \mathbb{R}_0^d:\kappa <\|z\|\le 1\},&{\rm if}~\beta =1,\\
\{z\in \mathbb{R}_0^d:\|z\|>\kappa\},&{\rm if}~\beta\in (1,2),\\
\end{cases}
\]
we have by (\ref{long q}) that
\begin{eqnarray*}
\int_{\mathbb{R}_0^d}{\bf 1}_B(z)zh(T_{h^{-1/\beta}}\nu_{\sigma,q}^{\alpha,\beta})(dz)
&=&\int_{S^{d-1}}\xi\sigma (d\xi)\int_0^{\infty}{\bf
1}_B(h^{-1/\beta}r\xi)rh^{1-1/\beta} q(r,\xi)dr\\
&=&\int_{S^{d-1}}\xi\sigma (d\xi)\int_0^{\infty}{\bf
 1}_B(r\xi)rh^{1+1/\beta}q(h^{1/\beta}r,\xi)dr\\
&\to&\int_{S^{d-1}}\xi\sigma_1(d\xi)\int_0^{\infty}{\bf 1}_B(r\xi)r
\frac{dr}{r^{\beta+1}},
\end{eqnarray*}
as $h\to \infty$, where the convergence holds true since for $h\in
 (1,\infty),$ and with the help of (\ref{long q}), 
\[
 \left\|\int_{S^{d-1}}\xi\sigma (d\xi)\int_0^{\infty}{\bf 1}_B(r\xi)r
h^{1+1/\beta}q(h^{1/\beta}r,\xi)dr\right\|\asymp \left\|\int_{S^{d-1}}
\xi\sigma (d\xi)\int_0^{\infty}{\bf 1}_B(r\xi)rq(r,\xi)dr\right\|<\infty.
\]
Next, let $(\alpha,\beta)\in [1,2)\times (0,1)$.
Then, observe that for each $\kappa >0$, and as $h\to \infty$,
\[
 -b_{\alpha,\beta}-\int_{\kappa<\|z\|\le 1}zh(T_{h^{-1/\alpha}}\nu_{\sigma,q}^{\alpha,\beta})(dz)\to-\int_{\|z\|>\kappa}z\nu_{\sigma_1}^{\alpha}(dz),
\]
where the convergence holds true as before.
This completes the proof of (i).

(ii) The random vector $h^{-1/2}X_{h}^{LS}$ is infinitely divisible
 with generating triplet
\[
 \left(-\int_{\|z\|\ge 1}z
 h(T_{h^{-1/2}}\nu^{\alpha,\beta}_{\sigma,q})(dz),0,h(T_{h^{-1/2}}\nu_{\sigma,q}^{\alpha,\beta})\right).
\]
Letting $f$ be a bounded continuous function from $\mathbb{R}_0^d$ to
 $\mathbb{R}$ such that $|f|\le C<\infty$ and $f(z)\equiv 0$ on $\{z\in
 \mathbb{R}^d:\|z\|\le \epsilon\}$, for some $\epsilon >0$,
the L\'evy measure $h(T_{h^{-1/2}}\nu_{\sigma,q}^{\alpha,\beta})$
 converges vaguely to zero as $h\to \infty$ since for sufficiently large
 $h>0$,
\begin{eqnarray}
\nonumber \left|\int_{\mathbb{R}_0^d}f(z)h(T_{h^{-1/2}}\nu^{\alpha
	   ,\beta}_{\sigma,q})(dz)\right|&=& \left|\int_{S^{d-1}}\sigma
 (d\xi)\int_{h^{1/2}\epsilon}^{\infty}f(h^{-1/2}r\xi)hq(r,\xi)dr\right|\\
\nonumber &\le&hC\int_{S^{d-1}}\sigma(d\xi)\int_{h^{1/2}\epsilon}^{\infty}
hq(r,\xi)dr\\
\nonumber &\asymp&hC\int_{S^{d-1}}c_2(\xi)\sigma(d\xi)\int_{h^{1/2}\epsilon}^{\infty}
\frac{dr}{r^{\beta+1}}\\
\label{long levy measure convergence} &=&h^{1-\beta/2}C\sigma_2(S^{d-1})\frac{\epsilon^{-\beta}}{\beta}\to 0,
\end{eqnarray}
as $h\to \infty$.
For the convergence of the Gaussian component, we have as $h\to +\infty$
 and for each $\kappa >0$,
\begin{equation}\label{gaussian long equation}
\int_{\|z\|\le
 \kappa}zz'h(T_{h^{-1/2}}\nu^{\alpha,\beta}_{\sigma,q})(dz)=
\int_{\|z\|\le h^{1/2}\kappa}zz'\nu_{\sigma,q}^{\alpha,\beta}(dz)\to \int_{\mathbb{R}_0^d}zz'\nu^{\alpha,\beta}_{\sigma,q}(dz),
\end{equation}
which is clearly well defined since $\int_{\mathbb{R}_0^d}\|z\|^2\nu^{\alpha,\beta}_{\sigma,q}(dz)<\infty$.
Finally, for sufficiently large $h>0$,
\begin{eqnarray}\nonumber
\left\|\int_{\|z\|>\kappa}zh(T_{h^{-1/2}}\nu^{\alpha,\beta}_{\sigma,q})
(dz)\right\|&=&h^{1/2}\left\|\int_{S^{d-1}}\xi\sigma
(d\xi)\int_{h^{1/2}\kappa}^{\infty}rq(r,\xi)dr\right\|\\
\nonumber &\asymp&h^{1/2}\int_{h^{1/2}\kappa}^{\infty}r\frac{dr}
{r^{\beta+1}}\left\|\int_{S^{d-1}}\xi\sigma_2(d\xi)\right\|\\
\label{long shift convergence}&=&h^{1-\beta/2}\frac{\kappa^{1-\beta}}
{\beta-1}\left\|\int_{S^{d-1}}\xi\sigma_2(d\xi)\right\|,
\end{eqnarray}
As $h\to \infty,$ (\ref{long shift convergence}) converges to zero and
 this concludes the proof of (ii).
\end{proof}

For $\beta =2$, layered stable processes do not seem to possess
 any nice long time behavior, and this can be seen from the
 improper convergence of the L\'evy measure, i.e., as $h\to \infty$,
 $h(T_{h^{-1/2}}\nu^{\alpha ,2}_{\sigma,q})$ converges vaguely to
\[
 \int_{S^{d-1}}\sigma_2(d\xi)\int_0^{\infty}{\bf 1}_B(r\xi)\frac{dr}{r^{2+1}},\quad
 B\in \mathcal{B}(\mathbb{R}_0^d),
\] 
which is not well defined as a L\'evy measure.
However, additional assumptions on $\sigma_2$ lead to the weak
convergence towards a Brownian motion as $\beta$ approaches to 2.

\begin{Prop}\label{special Gaussian convergence}
Let $\{X^{LS}_t:t\ge 0\}\sim LS_{\alpha,\beta}(\sigma ,q;0)$ in $\mathbb{R}^d$.

\noindent (i) Let $\beta \in (1,2)$ and let $\eta=-\int_{S^{d-1}}\xi
 \sigma (d\xi)\int_1^{\infty}rq(r,\xi)dr$.
If $\sigma_2$ is uniform on $S^{d-1}$ such that
 $\sigma_2(S^{d-1})=d(2-\beta)$, then 
\[
 \{h^{-1/\beta}(X^{LS}_{ht}+ht\eta):t\ge 0\}\stackrel{d}{\to}
\{W_t:t\ge 0\},\quad as~h\to \infty,~\beta \uparrow 2,
\]
where $\{W_t:t\ge 0\}$ is a $d$-dimensional (centered) standard Brownian
 motion.
(The limit is taken over $h\to \infty$ first.)

\noindent (ii) Let $\beta \in (2,\infty )$ and let $\eta$ be the constant
 (\ref{def of constant a}). 
If $\sigma_2$ is symmetric such that $\sigma_2(S^{d-1})=\beta -2,$ then
\[
 \{h^{-1/2}(X_{ht}^{LS}+ht\eta):t\ge 0\}\stackrel{d}{\to}\{W_t:t\ge 0\},\quad
 as~h\to \infty,~\beta\downarrow 2,
\]
where $\{W_t:t\ge 0\}$ is a centered Brownian motion with covariance
 matrix $\int_{\mathbb{R}_0^d}zz'\nu_{\sigma,q}^{\alpha,2}(dz).$
(The limit can be taken either over $h\to \infty$ or over $\beta
 \downarrow 2$ first.)
\end{Prop}

\begin{proof}
(i) By Theorem \ref{short of eS} (i),
 $h^{-1/\beta}(X^{LS}_h+h\eta)\stackrel{\mathcal{L}}{\to}X^{(\beta)}_1$, as
 $h\to \infty$, where $\{X^{(\beta)}_t:t\ge 0\}\sim S_{\beta}(\sigma_2;0).$
Then, by E.18.7-18.8 of Sato \cite{sato}, we get
 $\mathbb{E}[e^{i\langle
 y,X^{(\beta)}_1\rangle}]=\exp[-c_{\beta,d}\|y\|^{\beta}],$ where 
\[
 c_{\beta,d}=\frac{\Gamma(d/2)\Gamma ((2-\beta)/2)}{2^{\beta}\beta\Gamma((\beta+d)/2)}\sigma_2(S^{d-1}).
\]
Taking $\beta\uparrow 2$ and since $\Gamma (x+1)=x\Gamma (x),$ $x>0$, we
 get the result.

(ii) In view of (\ref{long levy measure convergence}), we get
\[
 h^{1-\beta/2}C\sigma_2(S^{d-1})\frac{\epsilon^{-\beta}}{\beta}=
h^{1-\beta/2}C(\beta -2)\frac{\epsilon^{-\beta}}{\beta}\to 0,
\]
as $\beta \downarrow 2$, which shows that the L\'evy measure
 $h(T_{h^{-1/2}}\nu_{\sigma,q}^{\alpha,\beta})$ converges vaguely to zero.
Moreover, in view of (\ref{long shift convergence}), for sufficiently
 large $h>0$,
\[
 h^{1-\beta/2}\frac{\kappa^{1-\beta}}{\beta-1}\left\|\int_{S^{d-1}}\xi\sigma_2(d\xi)\right\|=0,
\]
by the symmetry of $\sigma_2.$
In view of (\ref{gaussian long equation}), it remains to show that
 $\int_{\mathbb{R}_0^d}\|z\|^2\nu_{\sigma,q}^{\alpha,2}(dz)<\infty.$
Observe that
\begin{eqnarray*}
\int_{\mathbb{R}_0^d}\|z\|^2\nu_{\sigma,q}^{\alpha,\beta}(dz)=\int_{S^{d-1}}\sigma(d\xi)\int_0^1r^2q(r,\xi)dr+\int_{S^{d-1}}\sigma(d\xi)\int_1^{\infty}r^2q(r,\xi)dr,
\end{eqnarray*}
and the first term of the right hand side above is clearly uniformly
 bounded in $\beta\in [2,\infty)$, and so is the second term, since for
 every $\beta \in [2,\infty),$
\[
 \int_{S^{d-1}}\sigma(d\xi)\int_1^{\infty}r^2q(r,\xi)dr\asymp
 \int_{S^{d-1}}c_2(\xi)\sigma (d\xi)\int_1^{\infty}r^2\frac{dr}{r^{\beta
 +1}}=\frac{\sigma_2(S^{d-1})}{\beta-2}=1.
\]
This concludes the proof in view of (\ref{eS long 2}).
\end{proof}

\begin{Rem}{\rm
The short time behavior (Theorem \ref{short of eS}) and the
 (non-Gaussian) long time behavior (Theorem \ref{lS long-time theorem}
 (i)) can also be inferred from the series representation (\ref{original
 series}).
For simplicity, consider the symmetric case.
Letting $X_t:=\sum_{i=1}^{\infty}\overleftarrow{q}(\Gamma_i/T,V_i)V_i{\bf
 1}(T_i\le t),$ we have
\[
 h^{-1/\alpha}X_{ht}=\sum_{i=1}^{\infty}h^{-1/\alpha}\overleftarrow{q}
(\Gamma_i/(hT),V_i)V_i{\bf 1}(hT_i\le ht),
\]
and so for each $u>0$ and each $\xi\in S^{d-1}$ such that
 $c_1(\xi)\in [0,\infty)$, bounded convergence gives
\begin{eqnarray*}
h^{-1/\alpha}\overleftarrow{q}(h^{-1}u,\xi)&=&h^{-1/\alpha}\inf
 \left\{r>0:\int_r^{\infty}q(s,\xi)ds<h^{-1}u\right\}\\
&=&\inf
 \left\{r>0:\int_r^{\infty}h^{1+1/\alpha}q(h^{1/\alpha}s,\xi)ds<u\right\}\\
&\to&\inf
 \left\{r>0:c_1(\xi)\int_r^{\infty}s^{-\alpha-1}ds<u\right\}=
\left(\frac{\alpha u}{c_1(\xi)}\right)^{-1/\alpha},
\end{eqnarray*}
as $h\to 0$, which is indeed an $\alpha$-stable shot noise.
The (non-Gaussian) long time behavior can be inferred just similarly. 
}\end{Rem}

\section{Absolute continuity with respect to short time limiting stable process}
Two L\'evy processes, which are mutually absolutely continuous, share any
almost sure local behavior.
The next theorem confirms this fact in relation with the short
time behavior result of Theorem \ref{short of eS}.
Indeed, given any layered stable process with respect to some
probability measure, one can find a probability measure under which the
layered stable process is identical in law to its short time limiting
stable process.
This result should be compared with Section 4 of Rosi\'nski
\cite{super rosinski}.

Recall that $c_1$ and $c_2$ are integrable (with respect
to $\sigma$) functions on $S^{d-1}$ appearing in (\ref{short q}) and
(\ref{long q}), while $\sigma_1$ and $\sigma_2$ are the finite positive
measures (\ref{def of sigma1}) and (\ref{def of sigma2}), respectively.
As before, we use the notation $\nu_{\sigma ,q}^{\alpha,\beta}$ for the
L\'evy measure of a layered stable process $X:=\{X_t:t\ge 0\}\sim
LS_{\alpha,\beta}(\sigma,q;\eta)$, while $\nu_{\sigma}^{\alpha}$ is the
measure (\ref{def of stable levy measure}). 

\begin{Thm}\label{absolute continuity theorem}
Let $\mathbb{P}$, $\mathbb{Q}$ and $\mathbb{T}$ be probability measures
 on $(\Omega ,\mathcal{F})$ such that under $\mathbb{P}$ the canonical
 process $\{X_t:t\ge 0\}$ is a L\'evy process in $\mathbb{R}^d$ with
 $\mathcal{L}(X_1)\sim LS_{\alpha,\beta}(\sigma,q;k_0)$, while under
 $\mathbb{Q}$ it is a L\'evy process with $\mathcal{L}(X_1)\sim
 S_{\alpha}(\sigma_1;k_1)$.
Moreover, when $\beta \in (0,2)$ and under $\mathbb{T}$, $\{X_t:t\ge
 0\}$ is a L\'evy process with $\mathcal{L}(X_1)\sim
 S_{\beta}(\sigma_2;\eta)$, for some $\eta\in \mathbb{R}^d$.
Then, 

\noindent (i) $\mathbb{P}|_{\mathcal{F}_t}$ and $\mathbb{Q}|_{\mathcal{F}_t}$ are
 mutually absolutely continuous for every $t >0$ if and only if
\[
 k_0-k_1=
\begin{cases}
\int_{S^{d-1}}\xi\sigma(d\xi)\int_0^1rq(r,\xi)dr,&\alpha \in (0,1),\\
\int_{S^{d-1}}\xi\sigma(d\xi)\int_0^1r(q(r,\xi)-c_1(\xi)r^{-\alpha-1})dr,&\alpha
 =1,\\
\frac{1}{\alpha
 -1}\int_{S^{d-1}}\xi \sigma_1(d\xi)+\int_{S^{d-1}}\xi\sigma(d\xi)\int_0^1r(q(r,\xi)-c_1(\xi)r^{-\alpha-1})dr,&\alpha \in (1,2).
\end{cases}
\]

\noindent (ii) If $\alpha \ne \beta$, then for any choice of $\eta\in \mathbb{R}^d$,
 $\mathbb{P}|_{\mathcal{F}_t}$ and $\mathbb{T}|_{\mathcal{F}_t}$ are
 singular for all $t>0$.

\noindent (iii) For each $t>0$, 
\[
 \frac{d\mathbb{Q}}{d\mathbb{P}}|_{\mathcal{F}_t}=e^{U_t},
\]
where $\{U_t:t\ge 0\}$ is a L\'evy process defined on $(\Omega
 ,\mathcal{F},\mathbb{P})$ by
\begin{eqnarray}
 U_t&:=&\lim_{\epsilon\downarrow 0}\sum_{\{s\in (0,t]:\|\Delta X_s\|>\epsilon\}}\Bigg[
\ln \left(\frac{q(\|\Delta X_s\|,\Delta X_s/\|\Delta X_s\|)}{c_1(\Delta
X_s/\|\Delta X_s\|)\|\Delta X_s\|^{-\alpha-1}}\right)\nonumber \\
&&\qquad \qquad \qquad \qquad -t (\nu_{\sigma,q}^{\alpha,\beta}-
\nu_{\sigma_1}^{\alpha})(\{z\in\mathbb{R}_0^d:\|z\|>\epsilon\})\Bigg].
\label{levy process u}
\end{eqnarray}
In the above right hand side, the convergence holds $\mathbb{P}$-a.s. uniformly
 in $t$ on every interval of positive length.
\end{Thm}

\begin{proof}
\noindent (i) By Theorem 33.1 and Remark 33.3 of Sato~\cite{sato}, it is necessary
 and sufficient to show that the following three conditions hold;
\begin{gather}
 \int_{\{z:|\varphi (z)|\le
 1\}}\varphi(z)^2\nu_{\sigma_1}^{\alpha}(dz)<\infty, \label{absolute
 continuity condition 1}\\
 \int_{\{z:\varphi (z)>1\}}e^{\varphi
 (z)}\nu_{\sigma_1}^{\alpha}(dz)<\infty, \label{absolute continuity
 condition 2}\\
 \int_{\{z:\varphi
 (z)<-1\}}\nu_{\sigma_1}^{\alpha}(dz)<\infty,\label{absolute
 continuity condition 3}
\end{gather}
where the function $\varphi:\mathbb{R}_0^d\to \mathbb{R}$ is defined by 
$(d\nu_{\sigma,q}^{\alpha
 ,\beta}/d\nu_{\sigma_1}^{\alpha})(z)=e^{\varphi (z)},$ that is,
\[
 \varphi (z)=\ln\left(\frac{q(\|z\|,z/\|z\|)}{c_1(z/\|z\|)\|z\|^{-\alpha-1}}
\right),\quad z\in \mathbb{R}_0^d.
\]
Now, observe that 
\begin{equation}\label{absolute continuity domain 1}
 \lim_{\|z\|\to 0}\varphi (z)=\lim_{\|z\|\to 0}\ln
  \left(\frac{c_1(z/\|z\|)\|z\|^{-\alpha-1}}
   {c_1(z/\|z\|)\|z\|^{-\alpha-1}}\right)=0,
\end{equation}
and that as $\|z\|\to\infty$,
\begin{equation}\label{absolute continuity domain 2}
 \varphi (z)\sim \ln \left(\frac{c_2(z/\|z\|)\|z\|^{-\beta-1}}
{c_1(z/\|z\|)\|z\|^{-\alpha-1}}\right)=
\ln \left(\frac{c_2(z/\|z\|)}{c_1(z/\|z\|)}\right)+(\alpha -\beta)\ln \|z\|\to
\begin{cases}
-\infty,&{\rm if}~\alpha <\beta,\\
+\infty,&{\rm if}~\alpha >\beta.
\end{cases}
\end{equation}
The conditions (\ref{absolute continuity condition 1}) and
 (\ref{absolute continuity condition 3}) are thus immediately
 satisfied, respectively, by (\ref{absolute continuity domain 1}) and
 (\ref{absolute continuity domain 2}) with $\alpha <\beta$.
In view of (\ref{absolute continuity domain 2}) with $\alpha >\beta$,
 the condition (\ref{absolute continuity condition 2}) is satisfied since
$\int_{\{z:\varphi (z)>1\}}e^{\varphi (z)}\nu_{\sigma_1}^{\alpha}(dz)$
 is bounded from above and below by constant multiples of $\int_{\|z\|>1}\frac{q(\|z\|,z/\|z\|)}{c_1(z/\|z\|)\|z\|^{-\alpha-1}}\nu_{\sigma_1}^{\alpha}(dz)
=\nu_{\sigma,q}^{\alpha,\beta}(\{z\in\mathbb{R}_0^d:\|z\|>1\})$.
When $\alpha=\beta \in (0,2)$, we have, by (\ref{absolute continuity
 domain 1}) and (\ref{absolute continuity domain 2}),
\[
\begin{cases}
\lim_{\|z\|\to 0}\varphi (z)=0,&\\
\lim_{\|z\|\to \infty}\varphi (z)=\lim_{\|z\|\to \infty}\ln
 \left(\frac{c_2(z/\|z\|)}{c_1(z/\|z\|)}\right)<\infty .&
\end{cases}
\]
The condition (\ref{absolute continuity condition 1}) is then satisfied
 since $\int_{\{z:|\varphi (z)|\le
 1\}}\varphi(z)^2\nu_{\sigma_1}^{\alpha}(dz)$ is bounded from above and
 below by constant multiples of
 $\int_{\|z\|>1}\varphi(z)^2\nu_{\sigma_1}^{\alpha}(dz)$, which is
 further bounded by $C \nu_{\sigma_1}^{\alpha}(\{z\in
 \mathbb{R}_0^d:\|z\|>1\})$ for some constant $C$.
The conditions (\ref{absolute continuity condition 2}) and
 (\ref{absolute continuity condition 3}) are also satisfied since the domains
 $\{z\in \mathbb{R}_0^d:\varphi (z)>1\}$ and $\{z\in
 \mathbb{R}_0^d:\varphi(z)<-1\}$ are contained in some compact sets of
 $\mathbb{R}_0^d$.

(ii) It suffices to show that either one of the following two
 conditions always fails;
\begin{gather}
 \int_{\{z:\psi (z)>1\}}e^{\psi (z)}\nu_{\sigma_2}^{\beta}(dz)<\infty,
 \label{absolute continuity condition 5}\\
 \int_{\{z:\psi (z)<-1\}}\nu_{\sigma_2}^{\beta}(dz)<\infty,\label{absolute
 continuity condition 6}
\end{gather}
where the function $\psi:S^{d-1}\to \mathbb{R}$ is defined by
 $(d\nu_{\sigma,q}^{\alpha ,\beta}/d\nu_{\sigma_2}^{\beta})(z)=e^{\psi
 (z)},$ that is,
\[
 \psi(z)=\ln \left(\frac{q(\|z\|,z/\|z\|)}
{c_2(z/\|z\|)\|z\|^{-\beta-1}}\right),\quad z\in \mathbb{R}_0^d.
\]
As in the proof of (i), observe that
\[
 \lim_{\|z\|\to \infty}\psi (z)=\lim_{\|z\|\to \infty}\ln \left(
 \frac{c_2(z/\|z\|)\|z\|^{-\alpha-1}}{c_2(z/\|z\|)\|z\|^{-\alpha-1}}\right)=0,
\]
and that as $\|z\|\to 0$,
\[
 \psi (z)\sim \ln \left(
 \frac{c_1(z/\|z\|)\|z\|^{-\alpha-1}}{c_2(z/\|z\|)\|z\|^{-\beta-1}}\right)=
\ln \left(\frac{c_1(z/\|z\|)}{c_2(z/\|z\|)}\right)+(\beta -\alpha)\ln \|z\|\to
\begin{cases}
+\infty,&{\rm if}~\alpha >\beta,\\
-\infty,&{\rm if}~\alpha <\beta.
\end{cases}
\]
Therefore, the condition (\ref{absolute continuity condition 5}) fails when
 $\alpha >\beta$ since
\[
 \int_{\{z:\psi (z)>1\}}e^{\psi
 (z)}\nu_{\sigma_2}^{\beta}(dz)=\nu_{\alpha,q}^{\alpha,\beta}
(\{z\in\mathbb{R}_0^d:\varphi (z)>1\})=+\infty,
\]
while (\ref{absolute continuity condition 6}) fails when $\alpha <\beta$ since
$\nu_{\sigma_2}^{\beta}(\{z\in\mathbb{R}_0^d:\psi (z)<-1\})=+\infty.$

(iii) This is a direct consequence of (i) with Theorem 33.2 of
 Sato~\cite{sato}.
\end{proof}

\begin{Rem}\label{absolute example}{\rm 
As in Remark \ref{special remark}, let
\[
 q(r,\xi)=\sigma(S^{d-1})^{-1}(r^{-\alpha-1}{\bf
 1}_{(0,1]}(r)+r^{-\beta-1}{\bf 1}_{(1,\infty)}(r)), ~\xi\in S^{d-1}.
\]
Then, the L\'evy process $\{U_t:t\ge 0\}$ given in (\ref{levy process
 u}) becomes
\[
 U_t=(\alpha-\beta)\sum_{\{s\in (0,t]:\|\Delta X_s\|>1\}}\ln (\|\Delta
 X_s\|)-t\left(\frac{1}{\beta}-\frac{1}{\alpha}\right)\sigma (S^{d-1}).
\]
Intuitively speaking, $(d\mathbb{Q}/d\mathbb{P})|_{\mathcal{F}_t}$ replaces all
 $\beta$-stable jumps of a layered stable process up to time $t$ (i.e., 
jumps with absolute size greater than $1$) by the corresponding
 $\alpha$-stable jumps without changing direction.
Moreover, when $\alpha <\beta$, the L\'evy measure $\nu$ of $\mathcal{L}(U_1)$
 is concentrated on $(-\infty,0)$ and is given by
\[
 \nu (-\infty,y)=\alpha^{-1}\sigma
 (S^{d-1})e^{\frac{\alpha}{\alpha-\beta}y},\quad y<0,
\]
while when $\alpha >\beta$, it is concentrated on $(0,\infty)$ and is
 given by
\[
 \nu (y,\infty)=\alpha^{-1}\sigma
 (S^{d-1})e^{\frac{\alpha}{\alpha-\beta}y},\quad y>0. 
\] 
 Let us next restate the absolute continuity result
 (Theorem \ref{absolute continuity theorem}) based on the fact that
 a series representation generates sample paths of a L\'evy process
 directly by generating every single jump.
For simplicity, we consider the symmetric case.
Let $\{Y_t:t\ge 0\}$ be an $\alpha$-stable process with
 $\mathcal{L}(Y_1)\sim S_{\alpha}(\sigma;k_1)$.
By Lemma \ref{alpha series}, there exists a version of
 $\{Y_t:t\in [0,T]\}$ given by
\[
 Y'_t=\sum_{i=1}^{\infty}\left(\frac{\alpha \Gamma_i}{\sigma
 (S^{d-1})T}\right)^{-1/\alpha}V_i{\bf 1}(T_i\le t)+k_1t.
\]
Also, let $\{X_t:t\ge 0\}$ be a layered stable process with
 $\mathcal{L}(X_1)\sim LS_{\alpha,\beta}(\sigma,q;k_0)$.
In view of the series representation (\ref{special series}), there
 exists a version of $\{X_t:t\in [0,T]\}$ given by
\begin{eqnarray*}
&& X'_t=\sum_{i=1}^{\infty}\Bigg[\left(\frac{\beta \Gamma_i}{\sigma
 (S^{d-1})T}\right)^{-1/\beta}{\bf 1}_{(0,\sigma
 (S^{d-1})T/\beta)}(\Gamma_i)\\
&&\qquad \qquad \qquad +\left(\frac{\alpha \Gamma_i}{\sigma
 (S^{d-1})T}+1-\frac{\alpha}{\beta}\right)^{-1/\alpha}{\bf 1}_{(\sigma
(S^{d-1})T/\beta,\infty)}(\Gamma_i)\Bigg]V_i{\bf 1}(T_i\le t)+k_0t,
\end{eqnarray*}
where all the random sequences are the same as those appearing in
 $\{Y'_t:t\in [0,T]\}$ above.
By Theorem \ref{absolute continuity theorem}, they are
 mutually absolutely continuous if and only if
\[
 k_0-k_1=
\begin{cases}
\frac{1}{\alpha-1}\int_{S^{d-1}}\xi\sigma_1(d\xi),&{\rm if}~\alpha \in
 (0,1)\cup (1,2),\\
0,&{\rm if}~\alpha =1.
\end{cases}
\]
We infer that the L\'evy process $\{U_t:t\in [0,T]\}$ in the
 Radon-Nykodym derivative of Theorem \ref{absolute continuity theorem}
 (iii), that is,
\[
 \frac{d\mathbb{Q}}{d\mathbb{P}}|_{\mathcal{F}_t}=e^{U_t},
\]
has a version given by
\[
 U'_t=-\frac{\alpha -\beta}{\alpha}\sum_{i=1}^{\infty}\ln
 \left(\frac{\alpha \Gamma_i}{\sigma (S^{d-1})T}\right){\bf 1}_{(0,\sigma
 (S^{d-1})T/\alpha]}(\Gamma_i){\bf 1}(T_i\le
 t)-t\left(\frac{1}{\beta}-\frac{1}{\alpha}\right)\sigma (S^{d-1}).
\]
As a direct consequence, we have
\[
 \mathbb{P}(X\in B)=\mathbb{E}_{\mathbb{P}}[e^{U'_T}{\bf 1}_B(Y')], \quad B\in \mathcal{B}(\mathbb{D}([0,T],\mathbb{R}^d)).
\]
Moreover, in view of Theorem 33.2 of Sato~\cite{sato},
\[
 \frac{d\mathbb{P}}{d\mathbb{Q}}|_{\mathcal{F}_t}=e^{-U_t},
\]
and so we can derive a version of $\{U_t:t\in [0,T]\}$ in terms of the
 jumps of the layered stable process as follows;
\[
 U''_t=-\frac{\alpha -\beta}{\beta}\sum_{i=1}^{\infty}\ln
 \left(\frac{\beta \Gamma_i}{\sigma (S^{d-1})T}\right){\bf 1}_{(0,\sigma
 (S^{d-1})T/\beta]}(\Gamma_i){\bf 1}(T_i\le
 t)-t\left(\frac{1}{\beta}-\frac{1}{\alpha}\right)\sigma (S^{d-1}).
\]
Similarly, we have 
\[
 \mathbb{Q}(Y\in B)=\mathbb{E}_{\mathbb{Q}}[e^{-U''_T}{\bf 1}_B(X')], \quad B\in \mathcal{B}(\mathbb{D}([0,T],\mathbb{R}^d)).
\]
}\end{Rem}

\section{Concluding remarks}

\noindent $\bullet$ The weak convergence towards a Brownian motion,
proved in Proposition \ref{special Gaussian convergence} (i), is
interesting in  the sense that a stable process with uniformly dependent
components converges in law to standard Brownian motion, i.e., with
independent components.
It is also interesting to see how a stable process with independent
components converges towards a Brownian motion.  
To this end, for $i=1,\ldots,d$, let $a_{i}\in [0,\infty)$, let 
\[
 b_{i+}:=(0,\ldots,0,+1,0,\ldots,0),\quad b_{i-}:=(0,\ldots,0,-1,0,\ldots,0),
\]
where $+1$ and $-1$ are located in the $i$-th component, and finally set
\[
 \sigma
 (d\xi):=\sum_{i=1}^d\frac{2-\alpha}{2}a_i(\delta_{b_{i+}}(d\xi)+
\delta_{b_{i-}}(d\xi)),\quad \xi\in S^{d-1},
\]
where $\delta$ is the Dirac measure.
Clearly, $\sigma$ is a symmetric finite positive measure on $S^{d-1}$.
Also, let $\{X_t^{(\alpha)}:t\ge 0\}\sim S_{\alpha}(\sigma;0).$
Then, if $y_i$ is the $i$-th component of $y$, we have, using
$\Gamma (x+1)=x\Gamma (x),$ $x>0,$ 
\begin{eqnarray*}
\mathbb{E}[e^{i\langle y,X_1^{(\alpha)}\rangle}]&=&
\exp\left[-\frac{\Gamma(1/2)\Gamma((2-\alpha)/2)}{2^{\alpha}\alpha\Gamma((1+\alpha)/2)}\int_{S^{d-1}}|\langle y,\xi\rangle
     |^{\alpha}\sigma(d\xi)\right]\\
&=&\exp\left[-\frac{1}{2}\sum_{i=1}^d\frac{\Gamma(1/2)\Gamma(1+(2-\alpha)/2)}{2^{\alpha-2}\alpha\Gamma((1+\alpha)/2)}a_i|y_i|^{\alpha}\right]\\
&\to&\exp\left[-\frac{1}{2}\sum_{i=1}^da_i|y_i|^2\right],\quad {\rm
 as}~\alpha \uparrow 2.
\end{eqnarray*}
Therefore, we get $\{X_t^{(\alpha)}:t\ge 0\}\stackrel{d}{\to}\{W_t:t\ge
0\}$ as $\alpha \uparrow 2$, where $\{W_t:t\ge 0\}$ is a Brownian motion
with covariance matrix
\[
\begin{pmatrix}
 a_1 & 0   & \ldots & 0 \\
           0 & a_2 & \ldots & 0 \\
      \vdots & \vdots & \ddots & \vdots \\
          0 & 0 & \ldots & a_d \\ 
\end{pmatrix}
.
\]

\vspace{1em}
\noindent $\bullet$ By making use of the absolute continuity of L\'evy measures, we can
 derive two more forms of series representations for a layered stable
 process induced by the L\'evy measure (\ref{original lS}), with $\alpha
 <\beta$.
With the notations of Theorem \ref{absolute continuity theorem}, we get
 for $z\in \mathbb{R}_0^d,$
\[
 \frac{d\nu_{\sigma,q}^{\alpha,\beta}}{d\nu_{\sigma}^{\alpha}}(z)={\bf
 1}_{(0,1]}(\|z\|)+\|z\|^{\alpha-\beta}{\bf 1}_{(1,\infty)}(\|z\|) \le 1,
\]
and
\[
 \frac{d\nu_{\sigma,q}^{\alpha,\beta}}{d\nu_{\sigma}^{\beta}}(z)=
\|z\|^{\beta-\alpha}{\bf 1}_{(0,1]}(\|z\|)+{\bf 1}_{(1,\infty)}(\|z\|)\le 1.
\]
Then, by the rejection method of Rosi\'nski \cite{rosinski2}, the summands
 $\{\overleftarrow{q}(\Gamma_i/T,V_i)V_i\}_{i\ge 1}$ in (\ref{special series})
can be respectively replaced by
\[
 \left\{\left(\frac{\alpha
 \Gamma_i}{\sigma(S^{d-1})T}\right)^{-1/\alpha}{\bf 1}\left(
\frac{d\nu_{\sigma,q}^{\alpha,\beta}}{d\nu_{\sigma}^{\alpha}}
\left(\left(\frac{\alpha\Gamma_i}{\sigma(S^{d-1})T}\right)^{-1/\alpha}
V_i\right)\ge U_i\right)V_i\right\}_{i\ge 1},
\]
and
\[
 \left\{\left(\frac{\beta \Gamma_i}{\sigma(S^{d-1})T}\right)^{-1/\beta}
{\bf 1}\left(\frac{d\nu_{\sigma,q}^{\alpha,\beta}}{d\nu_{\sigma}^{\beta}}
\left(\left(\frac{\beta\Gamma_i}{\sigma(S^{d-1})T}\right)^{-1/\beta}V_i\right)\ge
 U_i\right)V_i\right\}_{i\ge 1},
\]
where $\{U_i\}_{i\ge 1}$ is a sequence of iid uniform random variables
 on $[0,1]$, independent of all the other random sequences.

\vspace{1em}
\noindent $\bullet$ In similarity to the work presented in \cite{ftsm},
it is possible to define a notion of fractional layered stable motion (fLSm).
Then, as in \cite{ftsm}, fLSm will, in short time, be close to
fractional stable motion (with inner index $\alpha$) while in long time
it is close to either fractional Brownian motion (if $\beta >2$) or to
fractional stable motion (with index $\beta <2$).

\vspace{1em}
\noindent $\bullet$ Let us observe some sample paths of a layered
stable process, generated via the series representation (\ref{special
series}).
By Theorem \ref{short of eS} and \ref{lS long-time theorem}, the entire
situation is exhausted by the following three cases;
\begin{enumerate2}
\item[(i)] $\alpha <\beta <2,$

\item[(ii)] $\beta \in (2,\infty),$

\item[(iii)] $\alpha >\beta$ with $\beta \in (0,2)$. 
\end{enumerate2}

\begin{figure}
\begin{center}
\begin{tabular}{c}
      \resizebox{100mm}{58mm}{\includegraphics{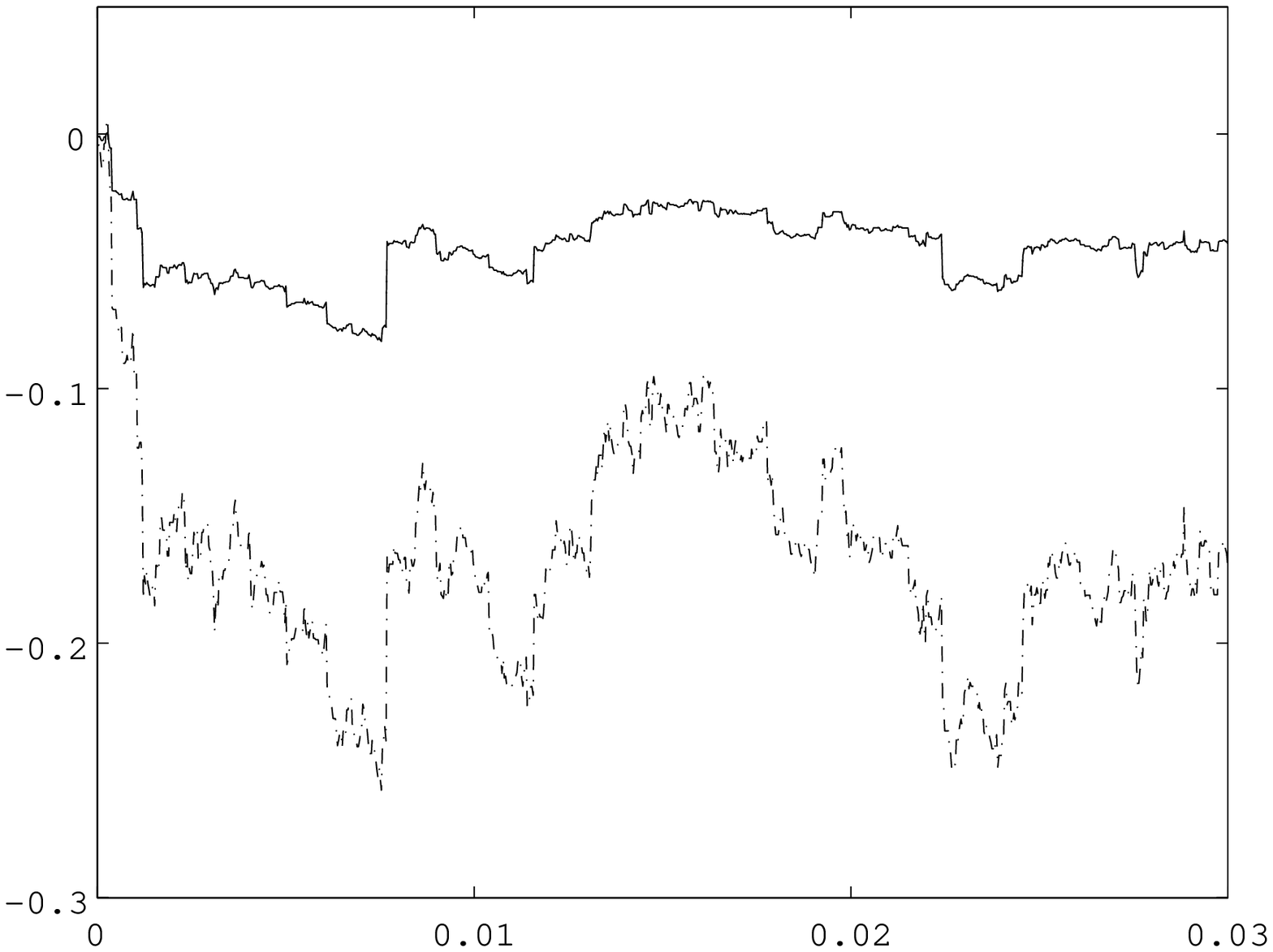}} \\
$t\in [0,0.03]$ \\
      \resizebox{100mm}{58mm}{\includegraphics{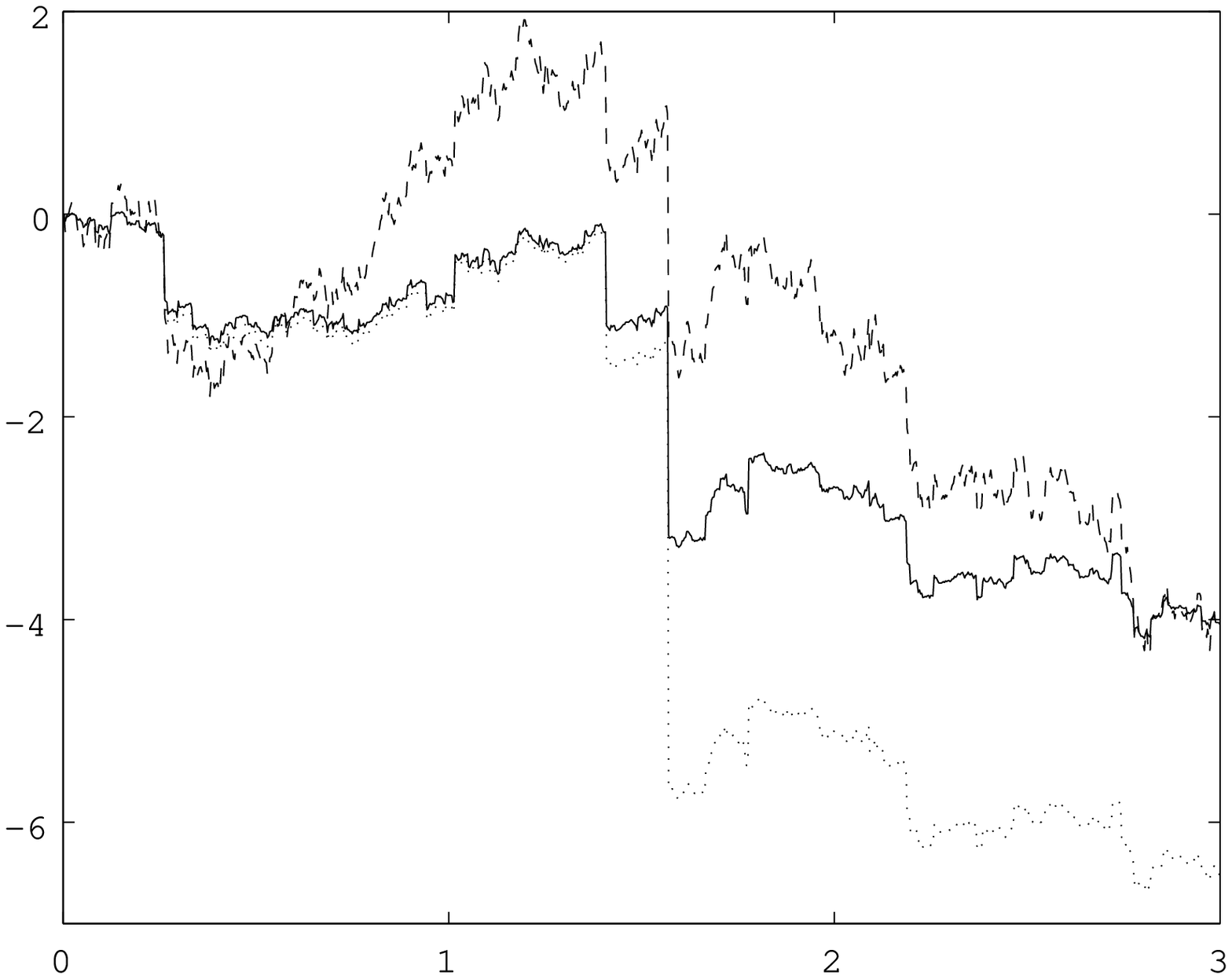}} \\
$t\in [0,3]$\\
      \resizebox{100mm}{58mm}{\includegraphics{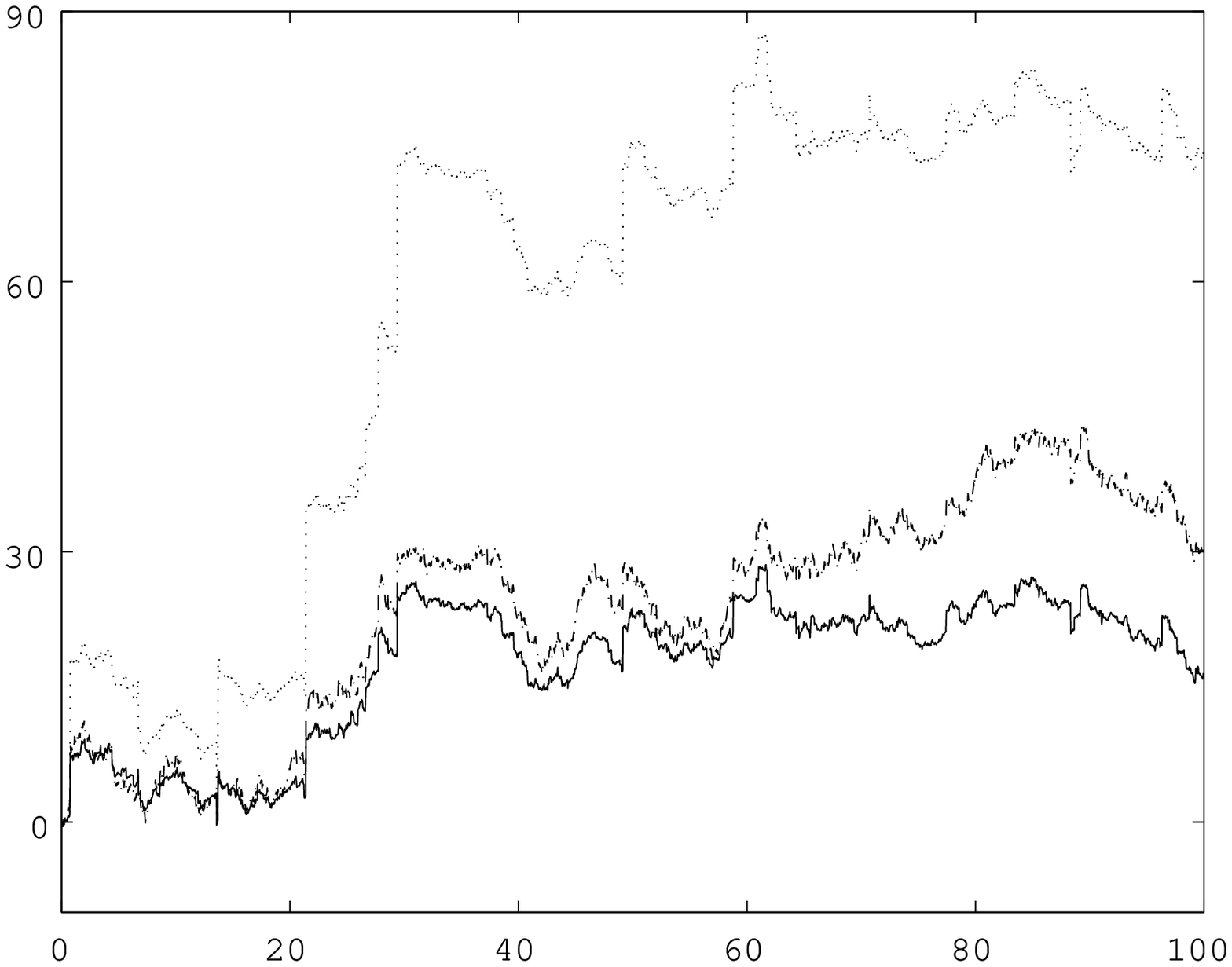}} \\
$t\in [0,100]$
\end{tabular}
\caption{Typical sample paths of layered stable process (---) with
 $(\alpha ,\beta )=(1.3,1.9)$, $1.3$-stable process ($\cdot \cdot \cdot$), and $1.9$-stable process (-$\cdot$-)}
\label{eS path}
\end{center}
\end{figure}

Figure \ref{eS path} corresponds to the case (i) and typical sample
paths of a symmetric layered stable process with
$(\alpha,\beta)=(1.3,1.9)$ are drawn in short, regular, and long time
span settings.
For better comparison, we also drew its corresponding $1.3$-stable and
$1.9$-stable processes.
All these sample paths are generated via the series representation
(\ref{special series}) for a layered stable process, or the one given in
Lemma \ref{alpha series} for stable processes.
Three sample paths within each figure are generated on a common
probability space in the sense that a common set of random sequences
$\{\Gamma_i\}_{i\ge 1}$, $\{V_i\}_{i\ge 1}$ and $\{T_i\}_{i\ge 1}$ are used.
The desired short and long time behaviors are apparent.
In the top figure, the layered stable process and its short time
limiting stable process are almost indistinguishable in a graphical
sense (of course, not in a probabilistic sense).

\begin{figure}
\begin{center}
\begin{tabular}{cc}
      \resizebox{100mm}{58mm}{\includegraphics{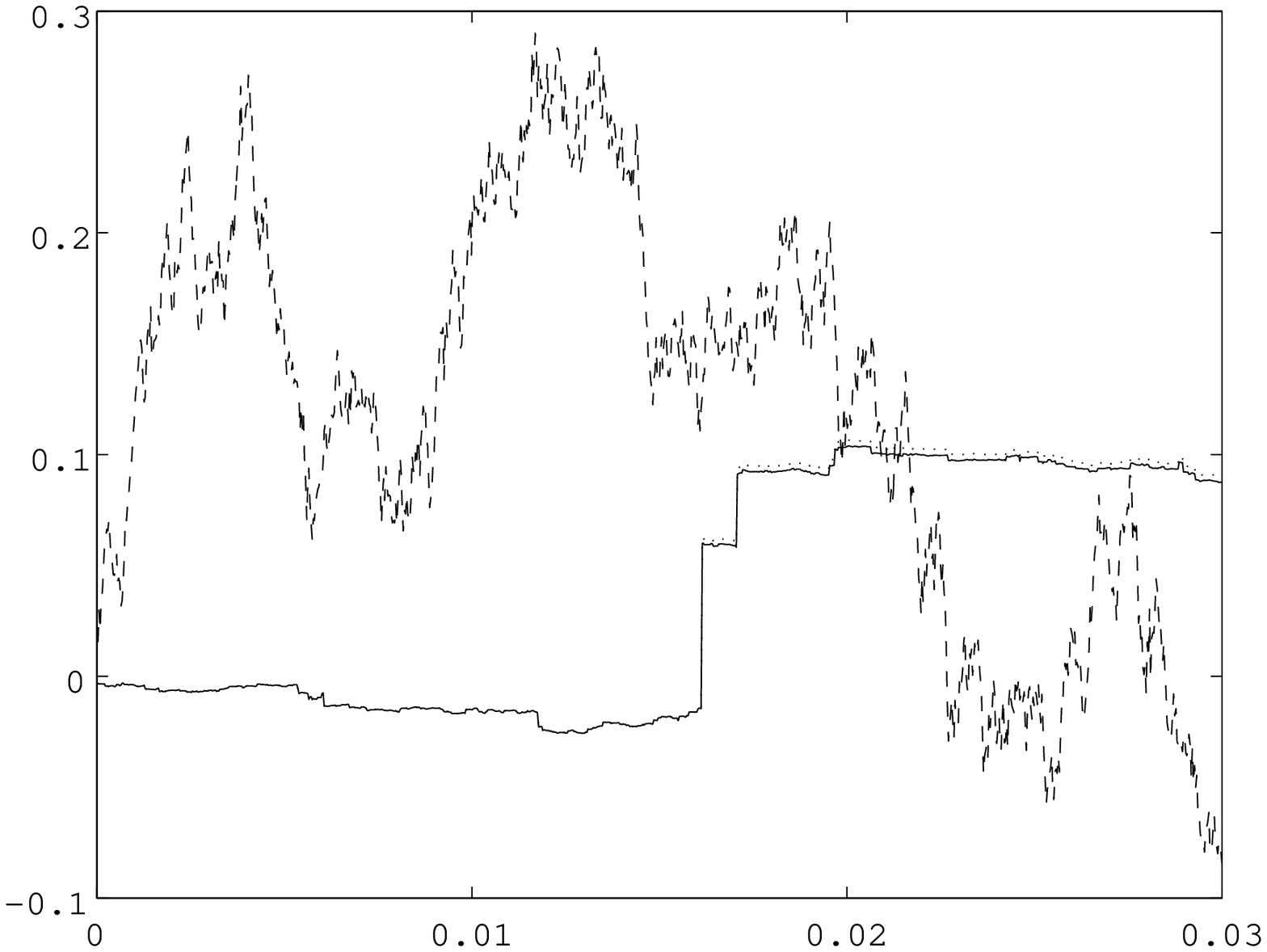}} \\
$t\in [0,0.03]$\\
      \resizebox{100mm}{58mm}{\includegraphics{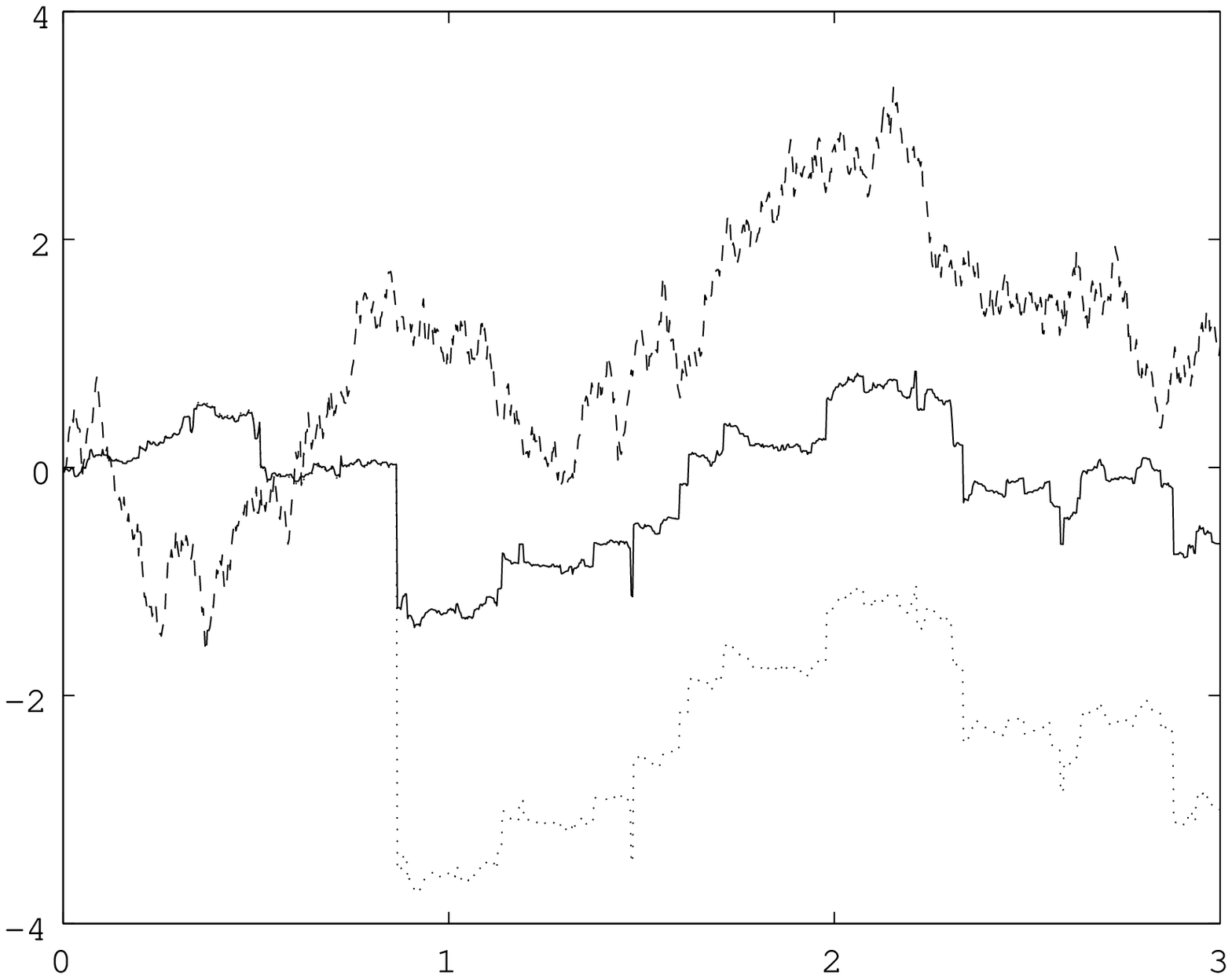}} \\
$t\in [0,3]$\\
      \resizebox{100mm}{58mm}{\includegraphics{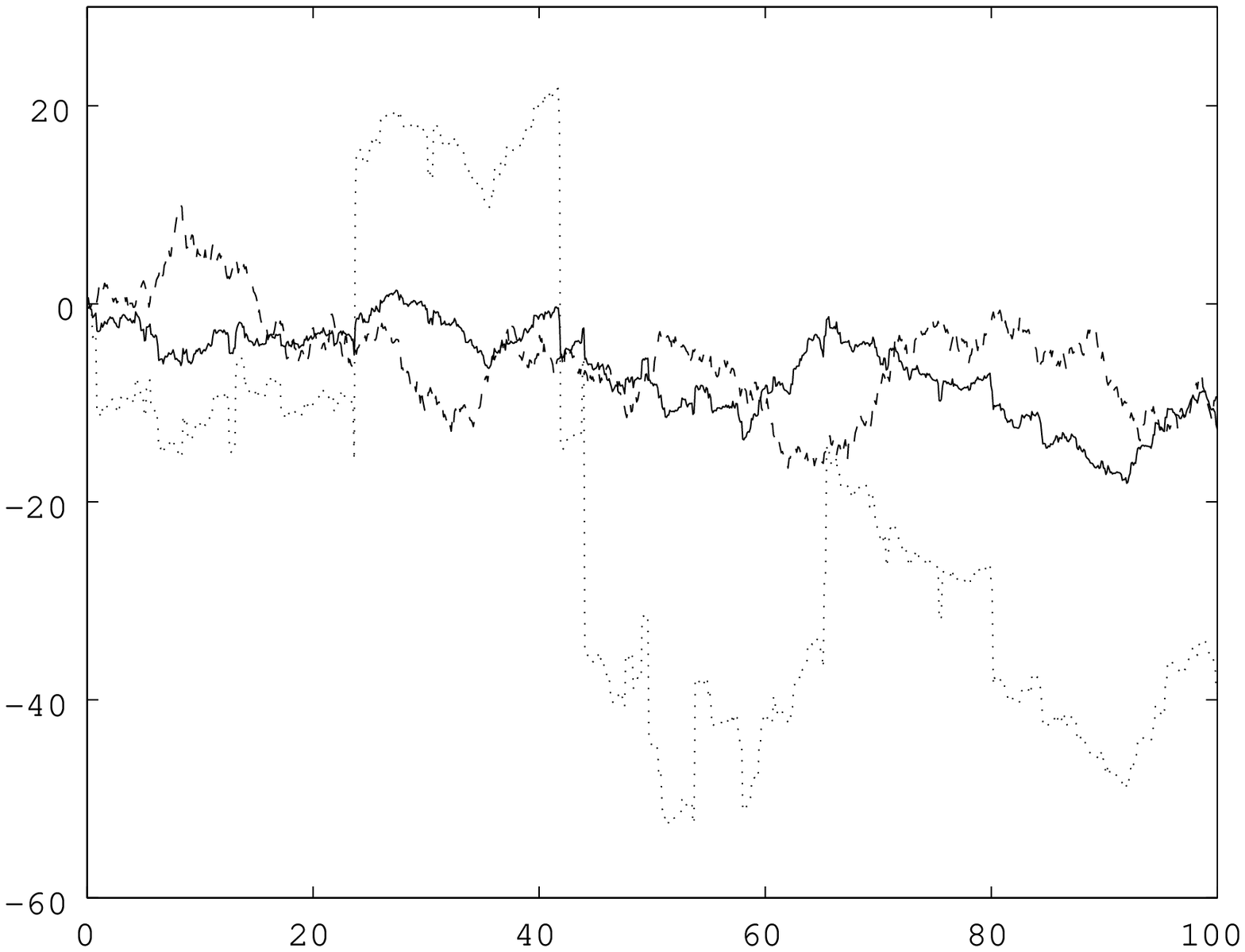}} \\
$t\in [0,100]$
\end{tabular}
\caption{Typical sample paths of layered stable process (---) with
 $(\alpha ,\beta )=(1.1,2.5)$, $1.1$-stable process ($\cdot \cdot
 \cdot$), and a Brownian motion (-$\cdot$-)}
\label{eS path 2}
\end{center}
\end{figure}

For the case (ii), we drew in Figure \ref{eS path 2}
typical sample  paths of a symmetric layered stable process with
$(\alpha,\beta)=(1.1,2.5)$, along with its corresponding $1.1$-stable
process and a Brownian motion with a suitable variance.
The layered stable process and the $1.1$-stable process are 
generated dependently as before, while the Brownian motion is
independent of the others.
As expected, the long time Gaussian type behavior (Theorem \ref{lS
long-time theorem} (ii)) is clearly apparent.
These stable type short time and Gaussian type long time behaviors have
long been considered to be very appealing in applications.
Such a study for asset price modeling will be presented elsewhere
\cite{finance paper}.

\begin{figure}
\begin{center}
\begin{tabular}{c}
      \resizebox{100mm}{58mm}{\includegraphics{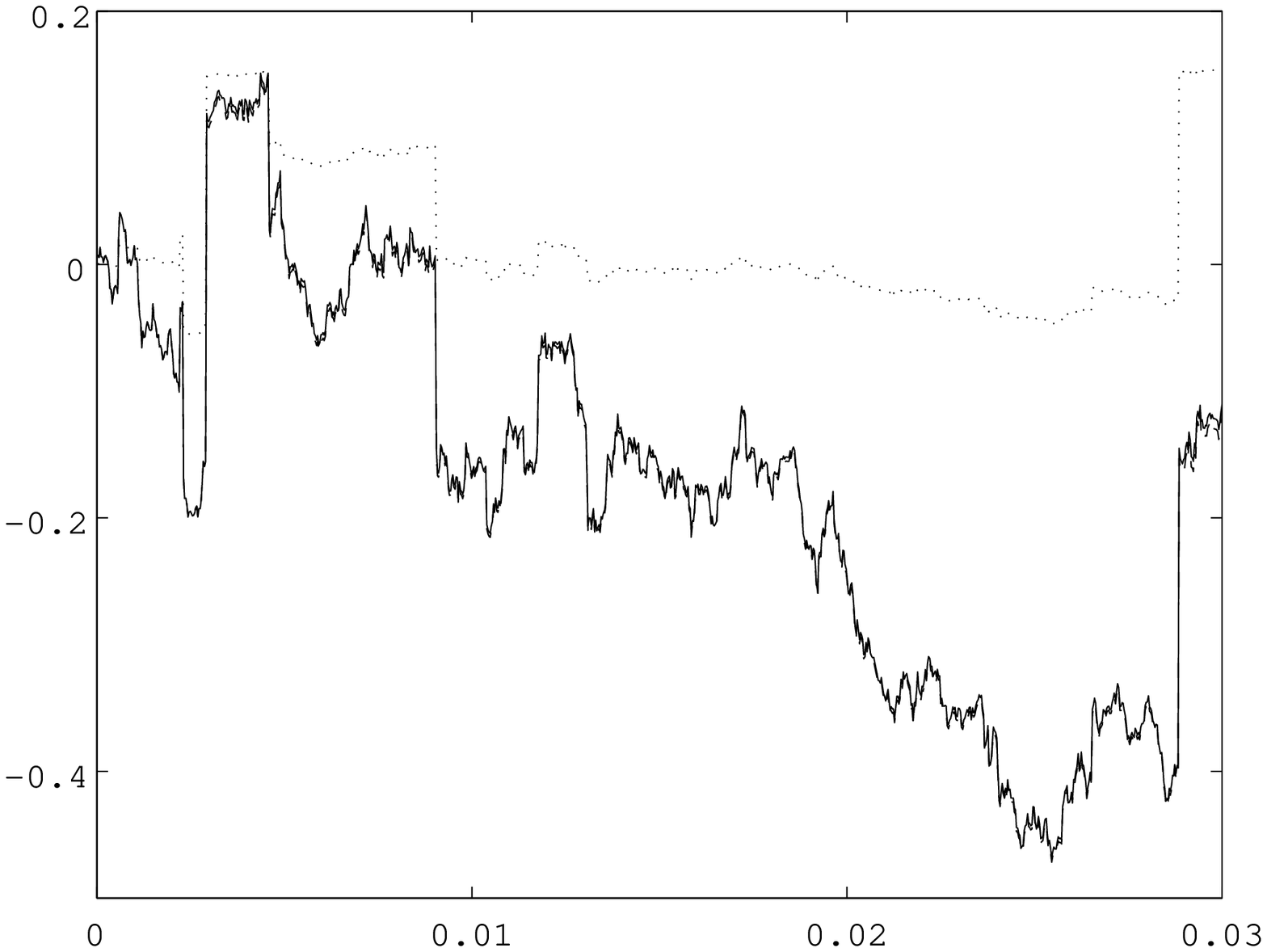}} \\
$t\in [0,0.03]$ \\
      \resizebox{100mm}{58mm}{\includegraphics{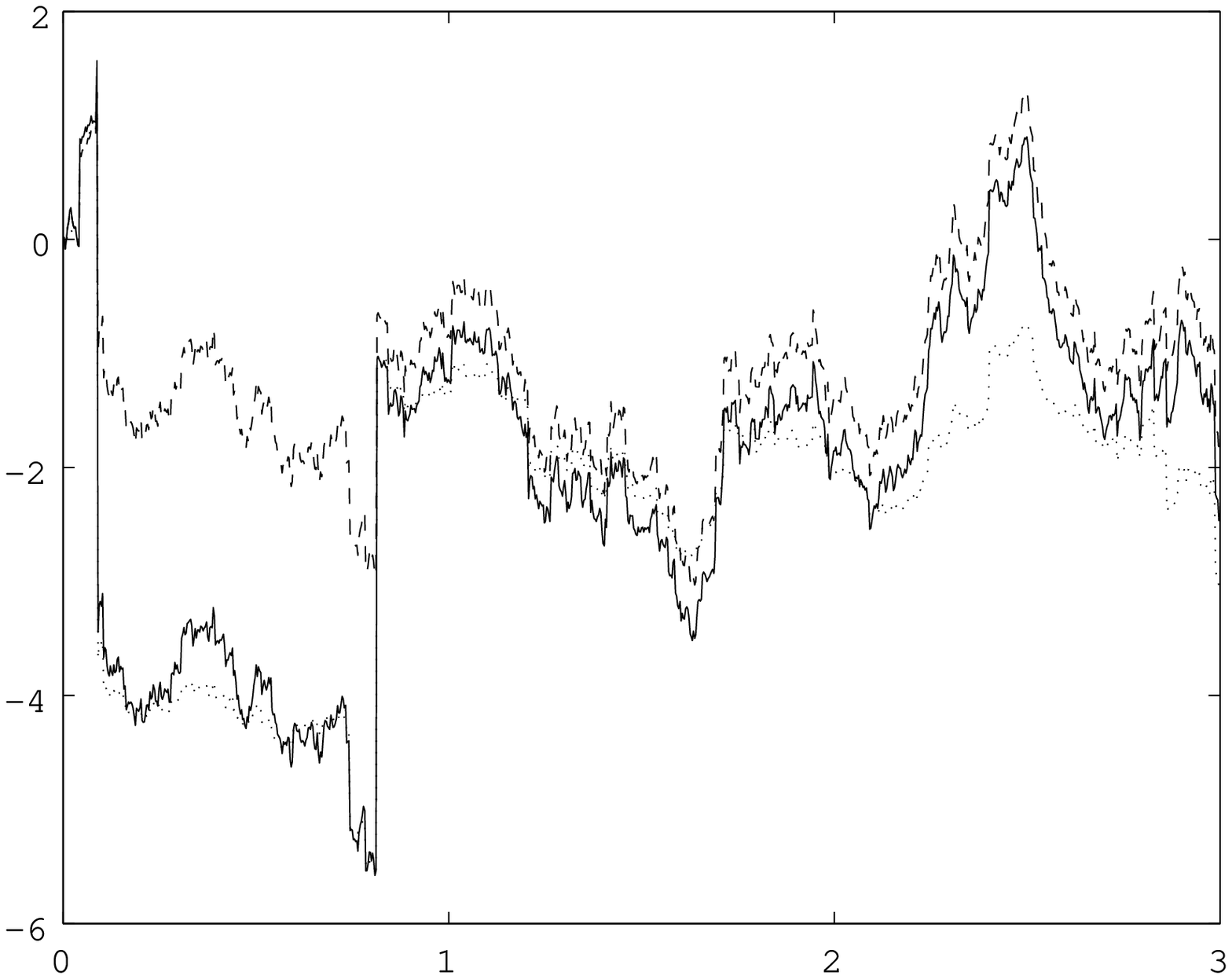}} \\
$t\in [0,3]$\\
      \resizebox{100mm}{58mm}{\includegraphics{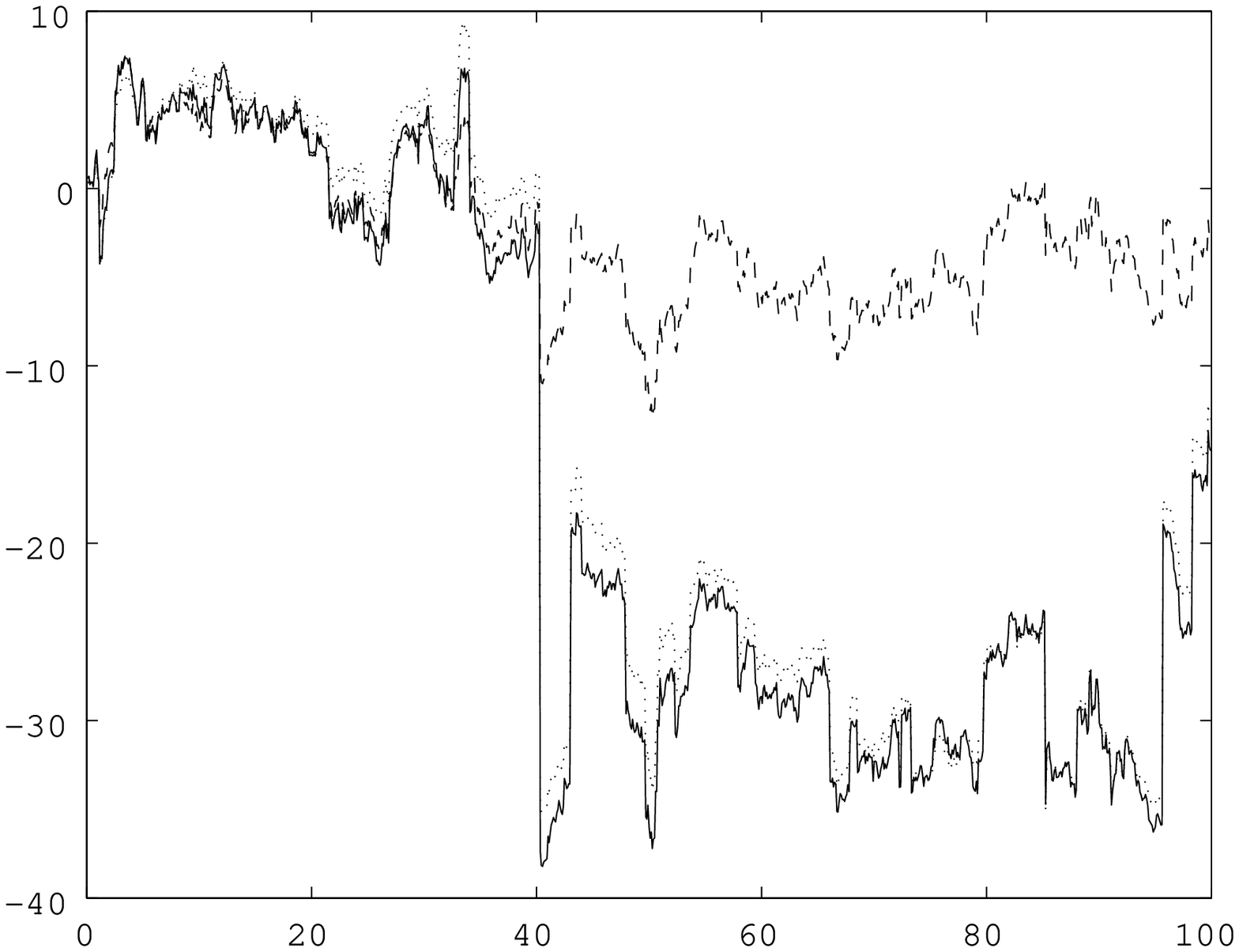}} \\
$t\in [0,100]$
\end{tabular}
\caption{Typical sample paths of layered stable process (---) with
 $(\alpha ,\beta )=(1.9,1.3)$, $1.9$-stable process (-$\cdot$-), and
 $1.3$-stable process ($\cdot \cdot \cdot$)}
\label{eS path 3}
\end{center}
\end{figure}

Finally, for the case (iii), we give in Figure \ref{eS path 3} typical
sample paths of a symmetric layered stable process with
$(\alpha,\beta)=(1.9,1.3)$, along with its corresponding $1.9$-stable
and $1.3$-stable processes.
Unlike the sample path behaviors observed in Figure \ref{eS path}, the
path of the layered stable processes behaves more continuously (like a
$1.9$-stable) in short time, while more discontinuously in long time (like
a $1.3$-stable).
In the short time figure, the layered stable and the $1.9$-stable are
graphically indistinguishable.

\vspace{1em}
\noindent $\bullet$ To finish this study, we briefly introduce another
generalization of stable processes.
Again, let $\mu$ be an infinitely divisible probability measure on
$\mathbb{R}^d$ and without Gaussian component.
Then, $\mu$ is {\it mixed stable} if its L\'evy measure is given by
\begin{equation}\label{mixed stable levy measure}
 \nu(B)=\int_{(0,2)}\int_{S^{d-1}}\sigma (d\xi)\int_0^{\infty}{\bf
 1}_B(r\xi)\frac{dr}{r^{\alpha+1}}\varphi (d\alpha),\quad B\in \mathcal{B}(\mathbb{R}_0^d),
\end{equation}
where $\varphi$ is a probability measure on $(0,2)$ such that
\[
 \int_{(0,2)}\frac{1}{\alpha (2-\alpha)}\varphi (d\alpha)<\infty.
\]
Its characteristic function is given by
\begin{eqnarray*}
\widehat{\mu}(y)&=&\exp\Bigg[i\langle y,\eta\rangle -\int_{(0,2)}c_{\alpha}\int_{S^{d-1}}|\langle
 y,\xi\rangle|^{\alpha}(1-i\tan \frac{\pi \alpha}{2}{\rm sgn}\langle
 y,\xi\rangle )\sigma (d\xi)\varphi (d\alpha)\\
&&\qquad \qquad \qquad \qquad \qquad -\varphi (\{1\})c_1\int_{S^{d-1}}(|\langle y,\xi\rangle
 |+i\frac{2}{\pi}\langle y,\xi\rangle \ln |\langle y,\xi\rangle| )\sigma (d\xi)\Bigg],
\end{eqnarray*}
for some $\eta \in \mathbb{R}^d$, and where $c_{\alpha}=|\Gamma
(-\alpha)\cos \frac{\pi \alpha}{2}|$ when $\alpha \ne 1$ while $c_1=\pi/2.$
Recall that in Example \ref{special remark} we defined the classes
$L_m$, $m=0,1,\ldots$ 
Let also $L_{\infty}:=\cap_{m=0}^{\infty}L_m$.
It is proved in Sato \cite{sato class L} that an infinitely divisible
probability measure without Gaussian component is in $L_{\infty}$ if and
only if its L\'evy measure has the form (\ref{mixed stable levy
measure}).
We have seen in Example \ref{special remark} that an infinitely
divisible probability measure $\mu$ is in $L_0$ if and only if the
L\'evy measure of $\mu$ has the form
\[
\int_{S^{d-1}}\sigma(d\xi)\int_0^{\infty}{\bf
 1}_B(r\xi)k_{\xi}(r)\frac{dr}{r},\quad B\in \mathcal{B}(\mathbb{R}_0^d),
\]
where $\sigma$ is a finite positive measure on $S^{d-1}$ and where
 $k_{\xi}(r)$ is a nonnegative function measurable in $\xi \in S^{d-1}$ and
 decreasing in $r>0.$   
Recently, Barndorff-Nielsen et al.\cite{bnms} defined the class
 $T$ by further requiring that the function $k_{\xi}(r)$ be completely
 monotone in $r$ for $\sigma$-a.e. $\xi$.
Mixed stable distributions are indeed in the class $T$ since
 $\int_{(0,2)}r^{-\alpha}\varphi (d\alpha)$ is completely monotone.

Finally, note that the associated L\'evy process that we call a {\it
mixed stable process} possesses an interesting series representation.  
For simplicity, assume that $\sigma$ in (\ref{mixed stable levy
measure}) is symmetric.
Let $\{\Gamma_i\}_{i\ge 1},$ $\{T_i\}_{i\ge 1}$ and $\{V_i\}_{i\ge 1}$
be random sequences defined as before.
In addition, let $\{\alpha_i\}_{i\ge 1}$ be a sequence of iid random
variables with common distribution $\varphi$.
Assume moreover that all these random sequences are mutually independent.
Then, with the help of the generalized shot noise method of Rosi\'nski
\cite{rosinski2}, it can be shown that the stochastic process
\[
 \left\{\sum_{i=1}^{\infty}\left(\frac{\alpha_i \Gamma_i}{\sigma
 (S^{d-1})T}\right)^{-1/\alpha_i}V_i{\bf 1}(T_i\le t):t\in [0,T]\right\},
\] 
converges almost surely uniformly in $t$ to a mixed stable process whose
marginal law at time $1$ is mixed stable with the L\'evy measure
(\ref{mixed stable levy measure}).
Comparing this result with the series representation of a stable process
given in Lemma \ref{alpha series}, a mixed stable process can be thought
of as a stable process with each of its jumps obeying a randomly chosen
stability index.

\end{document}